\documentclass{tac}
\usepackage{amsmath,amsfonts,graphics,ifpdf}

\author{George M. Bergman}

\thanks{%
arXiv: 1108.0958.
After publication of this note, updates, errata, related references
etc., if found, will be recorded at
\url{http://math.berkeley.edu/~gbergman/papers/}\,.
}

\address{%
Department of Mathematics,\\
University of California\\
Berkeley, CA 94720-3840, USA}

\title{On diagram-chasing in double complexes}

\copyrightyear{2012}


\eaddress{gbergman@math.berkeley.edu}

\ifpdf
  \usepackage{epic} % below: make dots same in epic, eepic
\newcommand{\middot}{\put(-1.0,0){\circle*{1.7}}}
\newcommand{\bigdot}{\circle*{2.6}}
  \newcommand{\diag}{\put(-2,0){\rotatebox{45}{\line(0,1){10}}}\put(2.5,3.5){\middot}\kern+.6em}
  \newcommand{\shortdiag}{\put(-5,0){\rotatebox{45}{\line(0,1){7}}}}
\else
  \usepackage{eepic}
  \newcommand{\middot}{\put(-6.0,-6.0){\mbox{\Huge$\cdot$}}}
  \newcommand{\bigdot}{\put(-3.6,-1.7){\mbox{\scriptsize$\bullet$}}}
  \newcommand{\diag}{\put(-2,7){\line(1,-1){7}}\put(2.0,3.5){\middot}\kern+.6em}
  \newcommand{\shortdiag}{\put(-5,5){\line(1,-1){5}}}
\fi %end of "ifpdf"

\newtheorem{theorem}{Theorem}
\newtheorem{lemma}[theorem]{Lemma}
\newtheorem{corollary}[theorem]{Corollary}

\newtheorem{definition}[theorem]{Definition}
\newtheorem{convention}[theorem]{Convention}

\newcommand{\cP}{{\scriptstyle\kern .20em\coprod\strut\kern -.08em}}

\newcommand{\bx}{{\kern.15em\line(1,0){3}\line(0,1){3}%
\raisebox{.3em}{\line(-1,0){3}\line(0,-1){3}\kern.3em}}}
\newcommand{\Bx}{{\kern.20em\line(1,0){4}\line(0,1){4}%
\raisebox{.4em}{\line(-1,0){4}\line(0,-1){4}\kern.4em}}}
\newcommand{\hor}{\put(0,3){\line(1,0){7}}\put(4.5,3.0){\middot}\kern+.7em}
\newcommand{\ver}{\put(2.0,-1.5){\line(0,1){9}}\put(3.0,3){\middot}\kern+.3em}
\newcommand{\A}{\mathbf{A}}
\newcommand{\B}{\mathbf{B}}
\newcommand{\Z}{\mathbb{Z}}
\newcommand{\N}{\mathbb{N}}
\newcommand{\bd}{\overline\delta}
\renewcommand{\r}{\mathrm}
\newcommand{\projLim}{\r{Lim}\raisebox{-.5em}{$\kern-1.8em{\longleftarrow}$}\,}
\newcommand{\injLim}{\r{Lim}\raisebox{-.5em}{$\kern-1.8em{\longrightarrow}$}\,}

\newcommand{\twosquare}[1]{
\begin{picture}(20,20)
\put(0,10){\line(0,1){10}}
\put(0,20){\line(1,0){10}}
\put(10,20){\line(0,-1){20}}
\put(10,0){\line(1,0){10}}
\put(20,0){\line(0,1){10}}
\put(20,10){\line(-1,0){20}}#1
\end{picture}}

\newcommand{\twocube}[1]{
\begin{picture}(30,30)
\multiput(15,0)(-15,15){2}{
\multiput(5,0)(-5,5){2}{
\multiput(0,0)(0,10){2}{\line(1,0){10}}
\multiput(0,0)(10,0){2}{\line(0,1){10}}
}
\multiput(5,0)(10,0){2}{\multiput(0,0)(0,10){2}{\shortdiag}}
}
#1
\end{picture}}

\begin{document}

\maketitle
\begin{abstract}
We construct, for any double complex in an abelian
category, certain ``short-distance''
maps, and an exact sequence involving these, instances of
which can be pieced together to give the
``long-distance'' maps and exact sequences
of results such as the Snake Lemma.

Further applications are given.
We also note what the building
blocks of an analogous study of triple complexes would be.
\end{abstract}
% - - - - - - - - - - - - - - - - - - - - - - - - - - - - - -

\section*{Introduction}
Diagram-chasing arguments frequently lead to ``magical''
relations between distant points of diagrams:  exactness implications,
connecting morphisms, etc..
These long connections are usually composites
of short ``unmagical'' connections, but the latter, and the
objects they join, are not visible in the proofs.
This note is aimed at remedying that situation.

Given a double complex in an abelian category, we will consider, for
each object $A$ of the complex, the familiar horizontal and vertical
homology objects at $A$ (which we will denote $A\hor$ and $A\ver\,),$
and two other objects, $A_\bx$ and $^\bx A$
(which we name the ``donor'' and the ``receptor'' at $A).$
For each arrow of the double complex, we prove
in \S\ref{S.salamander} the exactness of a certain $\!6\!$-term
sequence of maps between these objects (the ``Salamander Lemma'').
Standard results such as the $\!3\times 3$-Lemma, the Snake
Lemma, and the long exact sequence of homology associated with
a short exact sequence of complexes are recovered
in \S\S\ref{S.easy}-\ref{S.long} as easy
applications of this result.

In \S\ref{S.rows+cols} we generalize
the last of these examples, getting various exact diagrams
from double complexes with all but a few rows and columns exact.
The total homology of a double complex is examined
in \S\ref{S.total_homology}.

In \S\ref{SS.3-cx}
we take a brief look at the world of triple complexes,
and in \S\ref{SS.ratios} at
the relation between the methods of this note and J.\,Lambek's
homological formulation of Goursat's Lemma \cite{JL}.
We end with a couple of exercises.

\section{Definitions, and the Salamander Lemma}\label{S.salamander}
We shall work in an abelian category $\A.$
In the diagrams we draw, capital letters and boldface dots
\begin{picture}(8,8)
\put(4,4){\middot}
\end{picture}
will represent arbitrary objects of $\A.$
(Thus, such a dot does not imply a zero object,
but simply an object we do not name.)
Lower-case letters and arrows will denote morphisms.
When we give examples in categories of modules,
``module'' can mean left or right module, as the reader prefers.

A \emph{double complex} is an array of objects and
maps in $\A,$ of the form
\begin{equation}\begin{minipage}[c]{12pc}\label{d.2cx}
\begin{picture}(100,115)
\multiput(20,30)(0,20){4}{
\multiput(0,0)(20,0){4}{
\put(25,00){\middot}
\put(09,00){\vector(1,0){12}}
\put(24,16){\vector(0,-1){12}}}
\put(89,0){\vector(1,0){12}}
}
\multiput(20,10)(20,0){4}{
\put(24,16){\vector(0,-1){12}}}

\end{picture}\end{minipage}\end{equation}
extending infinitely on all sides, in which every row and every
column is a complex (i.e., successive arrows compose to zero), and all
squares commute.
Note that a ``partial'' double complex such as
\begin{equation}\begin{minipage}[c]{13pc}\label{d.4lem}
\begin{picture}(130,40)
\multiput(30,0)(30,0){3}{
\put(07,00){\middot}
\put(14,00){\vector(1,0){17}}
\put(6.5,24){\vector(0,-1){17}}
\put(07,30){\middot}
\put(14,30){\vector(1,0){17}}}
\put(127,00){\middot}
\put(126.5,24){\vector(0,-1){17}}
\put(127,30){\middot}
\end{picture}
\vspace{1em}
\end{minipage}\end{equation}
can be made a double complex by completing it with zeroes
on all sides; or by writing in some kernels and cokernels,
and then zeroes beyond these.
Thus, results on double complexes will be applicable to such finite
diagrams.

Topologists often prefer double complexes with \emph{anti}commuting
squares; but either sort of double complex can be turned
into the other by reversing the signs of the arrows in every other row.
In the theory of spectral sequences, vertical arrows generally go
upward, while in results like the Four Lemma they are generally
drawn downward; I shall follow the latter convention.

\begin{definition}\label{D.4things}
Let $A$ be an object of a double complex, with nearby
maps labeled as shown below:
\begin{equation}\begin{minipage}[c]{30pc}\label{d.4things}
\begin{picture}(90,80)
\put(10,10){
\put(07,00){\middot}
\put(14,00){\vector(1,0){17}}
\put(37,00){\middot}
\put(44,00){\vector(1,0){17}}\put(47.5,03){$h$}
\put(67,00){\middot}
\multiput(6.5,24)(30,0){3}{\vector(0,-1){17}}
\put(26,14){$f$}\put(57,14){$g$}
\put(07,30){\middot}
\put(14,30){\vector(1,0){17}}\put(17.5,33){$d$}
\put(31,27.5){$A$}
\put(44,30){\vector(1,0){17}}\put(48,33){$e$}
\put(67,30){\middot}
\multiput(6.5,54)(30,0){3}{\vector(0,-1){17}}
\put(-03,44){$b$}\put(27,44){$c$}
\put(07,60){\middot}
\put(14,60){\vector(1,0){17}}\put(18,63){$a$}
\put(37,60){\middot}
\put(44,60){\vector(1,0){17}}
\put(67,60){\middot}
\put(72,28){,}
\put(90,28){and let $p=ca=db,$ $q=ge=hf.$}
}
\end{picture}
\end{minipage}\end{equation}
Then we define\vspace{.2em}

$A\hor =\ \r{Ker}\ e\ /\ \r{Im}\ d,$ the horizontal homology
object at $A,$\vspace{.2em}

$A\ver\ =\ \r{Ker}\ f\ /\ \r{Im}\ c,$ the vertical homology
object at $A,$\vspace{.2em}

$^\bx A\ =\ (\r{Ker}\ e\,\cap\,\r{Ker}\ f)\,/\,\r{Im}\ p,$
which we shall call the \emph{receptor} at $A,$\vspace{.2em}

$A_\bx\ =\ \r{Ker}\ q\,/\,(\r{Im}\ c+\r{Im}\ d),$
which we shall call the \emph{donor} at $A.$
\end{definition}

From the inclusion relations among the kernels and images
in Definition~\ref{D.4things}, we get

\begin{lemma}\label{L.intermur}
For every object $A$ of a double complex, the identity map of $A$
induces a commuting diagram of maps:
\begin{equation}\begin{minipage}[c]{10pc}\label{d.intramural}
\begin{picture}(50,50)
\put(36,00){$A_\bx$}
\put(26,19){\vector(1,-1){12}}
\put(62,19){\vector(-1,-1){12}}
\put(14,23){$A\ver$}
\put(58,23){$A\hor$}
\put(38,44){\vector(-1,-1){12}}
\put(50,44){\vector(1,-1){12}}
\put(32,46){$^\bx A$}
\end{picture}
\endproof
\end{minipage}\end{equation}
\end{lemma}

\begin{definition}\label{D.intramural}
We shall call the maps shown in~\textup{(\ref{d.intramural})}
the \emph{intramural} maps associated with the object $A.$
\end{definition}

When we draw the diagram of a double complex, the donor and
receptor at an object will generally be indicated by small
squares to the lower right and upper left of the dot or letter
representing that object, as in~(\ref{d.extramural}) below.
Thus, the direction in which the square is displaced from
the letter is toward the most distant point of the diagram involved
in the definition of the object; in~(\ref{d.4things}), the domain or
codomain of the composite arrow $p$ or $q.$
(Of course, if one prefers to draw double complexes with arrows going
\emph{upward} and to the right, one should write $_\bx\kern.1em A$ and
$A^\bx$ for the receptor and donor at $A.)$

I will occasionally indicate horizontal or vertical homology
objects in a
diagram by marks $\hor$ and $\ver$ placed at the location of the
object; but this requires suppressing the symbol for the object itself.

The next result, whose proof is again straightforward,
motivates the names ``donor'' and ``receptor''.

\begin{lemma}\label{L.extramural}
Each arrow $f: A\to B$ in a double complex induces an
arrow $A_\bx\to{^\bx B}:$
\begin{equation}\begin{minipage}[c]{25pc}\label{d.extramural}
\begin{picture}(300,115)

% horizontal diagram:
\multiput(11,55)(50,0){3}{\vector(1,0){28}}
\put(44,52){$A$} \put(55,41){$\Bx$}
\put(94,52){$B$} \put(83,63){$\Bx$}
\dottedline{2}(65,46)(84,61)
 	\put(85,62){\vector(4,3){0}}
\multiput(50,44)(50,0){2}{
	\multiput(0,0)(0,50){2}{\vector(0,-1){28}}}

% vertical diagram:
\multiput(225,119)(0,-50){3}{\vector(0,-1){28}}
\put(219,77){$A$} \put(230,67){\Bx}
\put(219,27){$B$} \put(208,38){\Bx}
\dottedline{2}(231,65)(219,47)
	\put(217,44){\vector(-3,-4){0}}
\multiput(186,80)(0,-50){2}{
	\multiput(0,0)(50,0){2}{\vector(1,0){28}}}
\end{picture}
\endproof
\end{minipage}\end{equation}
\end{lemma}
\begin{definition}\label{D.extramural}
We shall call the morphism of Lemma~\ref{L.extramural}
the \emph{extramural} map associated with $f.$
\end{definition}

The global picture of the extramural maps in a double complex is
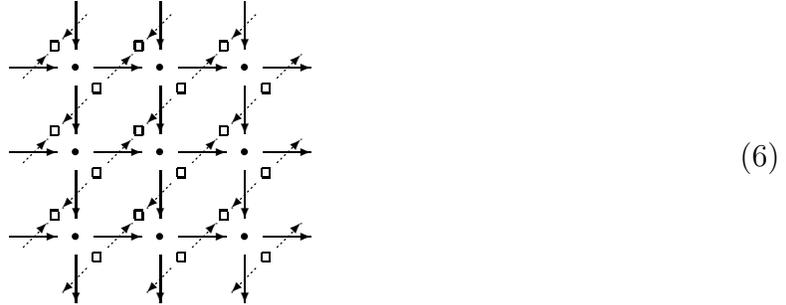
\begin{equation}\begin{minipage}[c]{35pc}\label{d.extra^n}
\begin{picture}(180,140)

\multiput(150,40)(0,32){3}{
\multiput(0,00)(32,0){4}{
\vector(1,0){18}
\dottedline{2}(-13,-4)(-5,04)
 	\put(-4,05){\vector(1,1){0}}
}
\multiput(28,00)(32,0){3}{
\put(-10.5,6.5){\bx}
\put(5.5,-9.5){\bx}
\put(1.5,0){\middot}
}}

\multiput(175,15)(32,0){3}{
\multiput(0,18)(0,32){4}{
\vector(0,-1){18}
\dottedline{2}(4,-5)(-4,-13)
 	\put(-5,-14){\vector(-1,-1){0}}
}}
\end{picture}\end{minipage}\end{equation}

In most of this note, I shall not use any notation
for these intramural and extramural maps.
Between any two of the objects we have constructed, we will not define
more than one map, so we shall be able to get by with an unlabeled
arrow, representing the unique map constructed between the
objects named.
(In \S\ref{S.total_homology}, where we will have more than one
map between the same pair of objects, I will introduce symbols
for some of these.)

To show a \emph{composite} of the maps we have defined, we may use
a long arrow marked with dots indicating the intermediate objects
involved, as in the statement of the following easily verified lemma.

\begin{lemma}\label{L.factor}
If $f:A\to B$ is a \emph{horizontal} arrow of a double complex
\textup{(}as in the first diagram of\textup{~(\ref{d.extramural}))},
then the natural induced map between
{\em vertical} homology objects, $A\ver\to B\ver,$ is the
composite of one extramural, and two intramural maps:
\begin{equation}\begin{minipage}[c]{35pc}\label{d.ver_to_ver}
\begin{picture}(400,25)
\put(103,10){$A\ver$}
\put(120,13){\vector(1,0){150}}
\put(157,02){$A_\bx$}
\put(165,13){\middot}
\put(207,18){$^\bx B$}
\put(215,13){\middot}
\put(272,10){$B\ver\,.$}
\end{picture}
\end{minipage}\end{equation}

Similarly, for a \emph{vertical} map $f:A\to B$ \textup{(}as in the
second diagram of\textup{~(\ref{d.extramural}))}, the
natural induced map of \emph{horizontal} homology objects is given by
\begin{equation}\begin{minipage}[c]{35pc}\label{d.hor_to_hor}
\begin{picture}(400,25)
\put(101,10){$A\hor$}
\put(122,13){\vector(1,0){148}}
\put(157,02){$A_\bx$}
\put(165,13){\middot}
\put(207,18){$^\bx B$}
\put(215,13){\middot}
\put(272,10){$B\hor\,.$}
\end{picture}
\endproof
\end{minipage}\end{equation}
\end{lemma}

We now come to our modest main result.
We will again state both the horizontal and vertical cases, since
we will have numerous occasions to use each.
The verifications are trivial if one is allowed to ``chase elements''.
To get the result in a general abelian category, one can use
one of the concretization theorems referred to
at \cite[Notes to Chapter VIII]{CW}, or
the method of generalized ``members'' developed in \cite[\S VIII.3]{CW},
or the related method discussed at \cite{elem_chase}.
(Regarding point~(vi) of \cite[\S VIII.4, Theorem~3]{CW}, we note that
this might be replaced by the
following more convenient statement, clear from the proof:
{\em Given $g:B\to C,$ and $x,y\in_m B$ with $gx\equiv gy,$
there exist $x'\equiv x,\ y'\equiv y$ with a common domain,
such that $gx'=gy'.)$}

We give our result a name analogous to that of the Snake Lemma.

\begin{equation}\begin{minipage}[c]{35pc}\label{d.salamander}
\begin{picture}(420,220)

% horizontal diagram:
\put(025,175){\vector(1,0){14}}
	\put(044,222){$C$}
	\put(050,214){\vector(0,-1){28}}
	\put(044,172){$A$}
	\put(050,164){\line(0,-1){14}}
\put(061,175){\vector(1,0){28}}
	\put(100,200){\line(0,-1){14}}
	\put(094,172){$B$}
	\put(100,164){\vector(0,-1){28}}
	\put(094,122){$D$}
\put(111,175){\line(1,0){14}}

% vertical diagram:
\put(070,100){\vector(0,-1){14}}
	\put(014,072){$C$}
	\put(031,075){\vector(1,0){28}}
	\put(064,072){$A$}
	\put(081,075){\line(1,0){14}}
\put(070,064){\vector(0,-1){28}}
	\put(045,025){\line(1,0){14}}
	\put(064,022){$B$}
	\put(081,025){\vector(1,0){28}}
	\put(114,022){$D$}
\put(070,014){\line(0,-1){14}}

% horizontal sketch:
	\put(350,225){\middot}
		\put(358,213){\bx}
	\dottedline{2}(350,214)(350,186)
		\put(338,183){\bx}
	\qbezier(357,213)(330,200)(346,177)
			\put(347,176){\vector(2,-3){0}}
	\put(345,172){$\hor$}
	\qbezier(351,170)(353,167)(357,165)
			\put(358,164){\vector(3,-2){0}}
		\put(358,163){\bx}
\dottedline{2}(361,175)(389,175)
	\qbezier(365,165)(371,168)(376,175)
	\qbezier(376,175)(381,182)(386,185)
			\put(387,185){\vector(2,1){0}}
		\put(388,183){\bx}
	\qbezier(395,185)(398,183)(400,180)
			\put(401,178){\vector(2,-3){0}}
	\put(395,172){$\hor$}
	\qbezier(406,173)(422,150)(395,137)
			\put(393,136){\vector(-2,-1){0}}
		\put(408,163){\bx}
	\dottedline{2}(400,164)(400,136)
	\put(400,125){\middot}
		\put(388,133){\bx}

% vertical sketch
	\put(320,075){\middot}
		\put(328,063){\bx}
	\dottedline{2}(331,075)(359,075)
		\put(358,083){\bx}
	\qbezier(332,068)(345,095)(368,079)
			\put(369,078){\vector(3,-2){0}}
	\put(367,072){$\ver$}
	\qbezier(374,074)(377,072)(379,068)
			\put(381,066){\vector(2,-3){0}}
		\put(378,063){\bx}
\dottedline{2}(370,064)(370,036)
	\qbezier(380,060)(377,054)(370,049)
	\qbezier(370,049)(363,044)(360,039)
			\put(359,037){\vector(-1,-2){0}}
		\put(358,033){\bx}
	\qbezier(361,030)(363,027)(366,025)
			\put(368,023){\vector(1,-1){0}}
	\put(367,022){$\ver$}
	\qbezier(372,019)(395,003)(408,030)
			\put(409,032){\vector(1,2){0}}
		\put(378,013){\bx}
	\dottedline{2}(381,025)(409,025)
	\put(420,025){\middot}
		\put(408,033){\bx}

% horizontal salamander sketch:
	\put(198,163){\middot}
		\put(206,151){\bx}
	\dottedline{2}(198,152)(198,124)
		\put(198,123){\vector(0,-1){0}}
		\put(186,121){\bx}
	\put(198,113){\middot} %
		\put(206,101){\bx}
\dottedline{2}(209,113)(237,113)
		\put(238,113){\vector(1,0){0}}
		\put(253,101){\bx}
	\put(248,113){\middot}
		\put(233,121){\bx}
	\dottedline{2}(248,102)(248,074)
		\put(248,073){\vector(0,-1){0}}
	\put(248,063){\middot}
		\put(236,071){\bx}
% the salamander itself:
	\qbezier(205,152)(183,131)(201,117)

 	\qbezier(201,117)(205,113)(210,113)
		\qbezier(210,113)(212,115)(217,120)
			\qbezier(213,116)(214,118)(216,121)
			\qbezier(213,116)(214,118.5)(215,122)
	\qbezier(210,113)(217,115)(223,119)
	\qbezier(223,119)(231,125)(238,125)
	\qbezier(238,125)(246,123)(253,119)
		\qbezier(258,114)(261,115)(261,125)
			\qbezier(261,120)(261.5,120)(263,125)
			\qbezier(261,120)(260.5,120)(259,125)
	\qbezier(253,119)(284,088)(235,067)
	\qbezier(235,067)(228,067)(225,071)
	\qbezier(225,071)(227,079)(230,079)
	\qbezier(230,079)(255,081)(249,101)
	\qbezier(249,101)(246,108)(238,108)
		\qbezier(246,105)(246,103)(241,097)
			\qbezier(244,101)(244,099)(243,096)
			\qbezier(244,101)(243,100)(239,097)
	\qbezier(238,108)(230,108)(223,103)
	\qbezier(223,103)(212,098)(206,099)
		\qbezier(207,099)(207,098)(199,094)
			\qbezier(203,096)(201,094)(199,092)
			\qbezier(203,096)(201,095)(198,096)
	\qbezier(206,099)(197,100)(192,108)
	\qbezier(192,108)(176,137)(205,152)
\end{picture}
\end{minipage}\end{equation}

\begin{lemma}[Salamander Lemma]\label{L.salamander}
Suppose $A\to B$ is a horizontal arrow in a double complex, and
$C,\ D$ are the objects above $A$ and below $B$ respectively,
as in the upper left diagram of\textup{~(\ref{d.salamander})}.
Then the following sequence
\textup{((\ref{d.salamander})}, upper right\textup{)},
formed from intramural and extramural maps, is exact:
\begin{equation}\begin{minipage}[c]{35pc}\label{d.sal_hor}
\begin{picture}(400,30)

\put(020,10){$C_\bx$}
\put(040,13){\vector(1,0){60}}
\put(063,18){$^\bx A$}
\put(065,13){\middot}
\put(105,10){$A\hor$}
\put(125,13){\vector(1,0){30}}
\put(160,10){$A_\bx$}
\put(180,13){\vector(1,0){30}}
\put(215,10){$^\bx B$}
\put(235,13){\vector(1,0){30}}
\put(270,10){$B\hor$}
\put(292,13){\vector(1,0){58}}
\put(314,02){$B_\bx$}
\put(316,13){\middot}
\put(355,10){$^\bx D.$}
\end{picture}
\end{minipage}\end{equation}

Likewise, if $A\to B$ is a vertical
arrow\textup{~((\ref{d.salamander})}, lower left\textup{)}, we have an
exact sequence \textup{((\ref{d.salamander})}, lower right\textup{)}:
\begin{equation}\begin{minipage}[c]{35pc}\label{d.sal_ver}
\begin{picture}(400,30)

\put(020,10){$C_\bx$}
\put(040,13){\vector(1,0){60}}
\put(063,18){$^\bx A$}
\put(065,13){\middot}
\put(105,10){$A\ver$}
\put(125,13){\vector(1,0){30}}
\put(160,10){$A_\bx$}
\put(180,13){\vector(1,0){30}}
\put(215,10){$^\bx B$}
\put(235,13){\vector(1,0){30}}
\put(270,10){$B\ver$}
\put(290,13){\vector(1,0){60}}
\put(313,02){$B_\bx$}
\put(315,13){\middot}
\put(355,10){$^\bx D.$}
\end{picture}
\end{minipage}\end{equation}
\endproof
\end{lemma}

In each case, we shall call the sequence displayed ``the $\!6\!$-term
exact sequence associated with the map $A\to B$ of the given
double complex''.

Remark: For $A\to B$ a horizontal arrow in a double complex,
and $C,$ $D$ as in the upper left diagram of~(\ref{d.salamander}),
not only~(\ref{d.sal_hor}) but also~(\ref{d.sal_ver}) makes sense,
but only the former is in general exact.
Indeed, by Lemma~\ref{L.factor}, the middle three maps
of~(\ref{d.sal_ver}) compose in that case to the natural
map from $A\ver$ to $B\ver,$ rather than to zero as they
would if it were an exact sequence.
Likewise, if $A\to B$ is a vertical map, then~(\ref{d.sal_hor})
and~(\ref{d.sal_ver}) both make sense, but only the
latter is in general exact.

\section{Special cases, and easy applications}\label{S.easy}

Note that the extramural arrows in~(\ref{d.extra^n}) stand
head-to-head and tail-to-tail, and so cannot be composed.
This difficulty is removed under appropriate conditions by the
following corollary to Lemma~\ref{L.salamander}, which one gets
on assuming the two homology objects in~(\ref{d.sal_hor})
or~(\ref{d.sal_ver}) to be $0:$

\begin{corollary}\label{C.extra-iso}
Let $A\to B$ be a horizontal \textup{(}respectively,
vertical\textup{)} arrow in a double complex, and suppose
the row \textup{(}resp., column\textup{)} containing this
map is exact at both $A$ and $B.$
Then the induced extramural map $A_\bx\to{^\bx B}$ is
an isomorphism.\endproof
\end{corollary}

Using another degenerate case of Lemma~\ref{L.salamander}, we
get conditions for \emph{intramural} maps to be isomorphisms, allowing
us to identify donor and receptor objects with classical
homology objects.

\begin{corollary}\label{C.inter-iso}
In each of the situations shown below, if the diagram is a
double complex, and the darkened row or column
\textup{(}the row or column through $B$ perpendicular to
the arrow connecting it with $A)$
is exact at $B,$ then the two intramural maps at $A$ indicated
above the diagram are isomorphisms.
\begin{equation}\begin{minipage}[c]{35pc}\label{d.inter-iso}
\begin{picture}(420,175)

\put(000,000){
\put(030,150){\put(-1,0){$^\bx A$}\put(17,003){\vector(1,0){16}}
	\put(17,005){$\cong$}\put(35,000){$A\hor$}}
\put(030,130){\put(0,0){$A\ver$}\put(17,003){\vector(1,0){16}}
	\put(17,005){$\cong$}\put(35,000){$A_\bx$}}
\put(035,034){\line(0,-1){11}}\put(075,032){\line(0,-1){11}}
\put(000,044){\put(-6,-3){$0$}\thicklines\put(5,0){\vector(1,0){20}}
	\put(28,-3){$B$}\put(47,0){\vector(1,0){20}}
	\put(75,0){\middot}\thinlines\put(85,0){\line(1,0){11}}}
\put(035,074){\vector(0,-1){20}} \put(075,075){\vector(0,-1){20}}
\put(000,084){\put(-6,-3){$0$}\put(5,0){\vector(1,0){20}}
	\put(28,-3){$A$}\put(47,0){\vector(1,0){20}}
	\put(75,0){\middot}\put(85,0){\line(1,0){11}}}
\put(035,103){\vector(0,-1){11}} \put(075,103){\vector(0,-1){11}}
}

\put(106,000){
\put(030,150){\put(-1,0){$^\bx A$}\put(17,003){\vector(1,0){16}}
	\put(17,005){$\cong$}\put(35,000){$A\ver$}}
\put(030,130){\put(0,0){$A\hor$}\put(20,003){\vector(1,0){14}}
	\put(18,005){$\cong$}\put(35,000){$A_\bx$}}
\put(035,034){\line(0,-1){11}}
	\put(075,034){\line(0,-1){11}}
\put(000,044){\put(16,0){\vector(1,0){11}}
	\put(35,0){\middot}\put(47,0){\vector(1,0){20}}
	\put(75,0){\middot}\put(85,0){\line(1,0){11}}}
\put(035,074){\vector(0,-1){20}}
	\thicklines\put(075,074){\vector(0,-1){20}}\thinlines
\put(000,084){\put(16,0){\vector(1,0){11}}
	\put(28,-3){$A$}\put(47,0){\vector(1,0){20}}
	\put(70,-3){$B$}\put(85,0){\line(1,0){11}}}
\put(035,108){\vector(0,-1){14}}
	\thicklines\put(075,108){\vector(0,-1){14}}\thinlines
\put(031,113){$0$} \put(071,113){$0$}
}

\put(207,000){
\put(030,150){\put(0,0){$A\hor$}\put(20,003){\vector(1,0){14}}
	\put(18,005){$\cong$}\put(36,000){$A_\bx$}}
\put(030,130){\put(-1,0){$^\bx A$}\put(17,003){\vector(1,0){16}}
	\put(17,005){$\cong$}\put(35,000){$A\ver$}}
\put(035,034){\line(0,-1){11}} \put(075,032){\line(0,-1){11}}
\put(000,044){\put(16,0){\vector(1,0){11}}
	\put(35,0){\middot}\put(47,0){\vector(1,0){20}}
	\put(69,-3){$A$}\put(85,0){\vector(1,0){20}}
	\put(110,-3){$0$}}
\put(035,074){\vector(0,-1){20}} \put(075,074){\vector(0,-1){20}}
\put(000,084){\put(16,0){\vector(1,0){11}}
	\put(37,0){\middot}\thicklines\put(47,0){\vector(1,0){20}}
	\put(69,-3){$B$}\put(85,0){\vector(1,0){20}}\thinlines
	\put(110,-3){$0$}}
\put(035,103){\vector(0,-1){11}} \put(075,103){\vector(0,-1){11}}
}

\put(333,000){
\put(030,150){\put(0,0){$A\ver$}\put(17,003){\vector(1,0){16}}
	\put(17,005){$\cong$}\put(35,000){$A_\bx$}}
\put(030,130){\put(-1,0){$^\bx A$}\put(17,003){\vector(1,0){16}}
	\put(17,005){$\cong$}\put(35,000){$A\hor$}}
\put(031,007){$0$} \put(071,007){$0$}
\put(035,034){\thicklines\vector(0,-1){14}}\thinlines
	\put(075,034){\vector(0,-1){14}}
\put(000,044){\put(16,0){\vector(1,0){11}}
	\put(28,-3){$B$}\put(47,0){\vector(1,0){20}}
	\put(70,-3){$A$}\put(85,0){\line(1,0){11}}}
\put(035,074){\thicklines\vector(0,-1){20}}\thinlines
	\put(075,074){\vector(0,-1){20}}
\put(000,084){\put(16,0){\vector(1,0){11}}
	\put(37,0){\middot}\put(47,0){\vector(1,0){20}}
	\put(75,0){\middot}\put(85,0){\line(1,0){11}}}
\put(035,103){\vector(0,-1){11}}\put(075,103){\vector(0,-1){11}}
}

\end{picture}
\end{minipage}\end{equation}
\end{corollary}

\proof
We will prove the isomorphism statements for the first
diagram; the proofs for the remaining diagrams are obtained
by reversing the roles of rows and columns, and/or reversing the
directions of all arrows.

In the first diagram, Corollary~\ref{C.extra-iso}, applied to
the arrow $0\to B,$ gives $^\bx B\cong 0_\bx = 0.$
The exact sequence~(\ref{d.sal_hor}) associated with the map
$0\to A$ therefore ends $0\to{^\bx A}\to A\hor\to 0,$
while the exact sequence~(\ref{d.sal_ver}) associated with the
map $A\to B$ begins $0\to A\ver\to A_\bx\to 0.$
This gives the two desired isomorphisms.
\endproof

{\em Remarks.}  Corollaries~\ref{C.extra-iso}
and~\ref{C.inter-iso} are easy to prove directly,
so in these two proofs, the Salamander Lemma was
not a necessary tool, but a guiding principle.

One should not be misled by apparent further dualizations.
For instance, one might think that if in the leftmost diagram
of~(\ref{d.inter-iso}), one assumed row-exactness at $A$ rather than
at $B,$ one would get isomorphisms $^\bx B\cong B\hor$
and $B\ver\cong B_\bx$ by reversing
the directions of vertical arrows, and applying the result proved.
However, there is no principle that allows us to reverse some
arrows (the vertical ones) without reversing others
(the horizontal arrows); and in fact, the asserted isomorphisms fail.
E.g., if the vertical arrow down from $B$ is the identity
map of a nonzero object,
and all other objects of the double complex are zero, then
row-exactness does hold at $A$ (since the row containing $A$ consists
of zero-objects), but the maps $^\bx B\to B\hor$
and $B\ver\to B_\bx$ both have the form $0\to B.$
(And if one conjectures instead
that one or both of the other two intramural
maps at $B$ should be isomorphisms,
the diagram with the horizontal map out of $B$ an identity map and
all other objects $0$ belies these statements.)
The common feature of the four situations of~(\ref{d.inter-iso}) that
is not shared by the results of reversing only vertical or
only horizontal arrows is that the arrow connecting
$A$ with $B,$ and the arrow connecting a zero object
with $B,$ have the same orientation relative to $B;$ i.e.,
either both go into it, or both come out of~it.

Most of the ``small'' diagram-chasing lemmas of homological
algebra can be obtained from the above two corollaries.
For example:
\begin{lemma}[The Sharp $3\times 3$ Lemma]\label{L.sharp}
In the diagram below \textup{(}\emph{excluding}, respectively
\emph{including} the parenthesized arrows; and ignoring for now the
boxes and dotted lines, which belong to the proof\textup{)},
if all columns, and all rows but the first,
are exact, then the first row \textup{(}again
excluding, respectively including the parenthesized arrow\textup{)}
is also exact.

\begin{equation}\begin{minipage}[c]{15pc}\label{d.sharp}
\begin{picture}(127,135)

\multiput(0,37)(0,32){3}{
\put(-3,-3){$0$}
\multiput(7,0)(32,0){3}{
\put(0,0){\vector(1,0){18}}
\put(25,25){\vector(0,-1){18}}}}
\multiput(28,128)(32,0){3}{$0$}

\multiput(107,69)(0,32){2}{
\put(-3,-3){$($}\put(3,0){\vector(1,0){15}}\put(18,-3){$0)$}}

\put(32,02){
\put(-6,-4){$\smile$}\put(-4,0){$0$}
\put(0,25){\vector(0,-1){15}}\put(-6,23){$\frown$}}

\put(26,34){
\put(00,00){$A''$}\put(00,32){$A$}\put(00,64){$A'$}
\put(32,00){$B''$}\put(32,32){$B$}\put(32,64){$B'$}
\put(64,00){$C''$}\put(64,32){$C$}\put(64,64){$C'$}}

\multiput(37,28)(0,32){3}{\bx}
\multiput(37,28)(16,16){4}{\dottedline{2}(6,6)(14,14)\put(16,16){\bx}}
\multiput(37,60)(16,16){2}{\dottedline{2}(6,6)(14,14)\put(16,16){\bx}}

\end{picture}\end{minipage}\end{equation}
\end{lemma}

\proof
We first note that in view of the exactness of the columns,
the first row of~(\ref{d.sharp}) consists of subobjects of
the objects of the second row, and restrictions of the maps
among these, hence it is, at least, a complex,
so the whole diagram is a double complex.

Corollaries~\ref{C.extra-iso} and~\ref{C.inter-iso},
combined with the exactness hypotheses
not involving the parenthesized arrows, now give us
\begin{equation}\begin{minipage}[c]{30pc}\label{d.AB_isos}
$A'\hor\,\cong\ A'_\bx\ \cong\ A'\ver\,=\ 0,\\[.3em]
B'\hor\,\cong\ B'_\bx\ \cong\ {^\bx B}
\ \cong\ A_\bx\ \cong\ A\ver\,=\ 0$
\quad(short dotted path in~(\ref{d.sharp})),
\end{minipage}\end{equation}
and, assuming also the exactness conditions involving the
parenthesized arrows,
\begin{equation}\begin{minipage}[c]{30pc}\label{d.C_isos}
$C'\hor\,\cong\ C'_\bx\ \cong\ ^\bx C\ \cong\ B_\bx
\ \cong\ ^\bx B''\ \cong \ A''_\bx\ \cong\ A''\ver\,=\ 0$\\
(long dotted path in~(\ref{d.sharp})).
\end{minipage}\end{equation}

(As an example of how to determine which statement of which of those
corollaries to use in each case, consider the first
isomorphism of the first line of~(\ref{d.AB_isos}).
Since it concerns an intramural isomorphism, it must be
an application of Corollary~\ref{C.inter-iso}.
Since the one of $^\bx A',$ $A'_\bx$ that it involves faces
away from the zeroes of the diagram, it must come from the
second row of isomorphisms in that corollary; and given
that fact, since it involves a horizontal homology object, it must
call on exactness at an object horizontally displaced from $A',$
hence must come from the second or fourth diagram
of~(\ref{d.inter-iso}).
Looking at the placement of the zeroes, we see that it must come from
the second diagram, and that the needed hypothesis, vertical
exactness at the object to the right of $A',$ is indeed present.
The isomorphisms other than the first and
the last in each line of~(\ref{d.AB_isos}) and in~(\ref{d.C_isos}),
corresponding to extramural maps,
follow from Corollary~\ref{C.extra-iso}.)

The consequent
triviality of the two (respectively three) horizontal homology
objects with which~(\ref{d.AB_isos}) (and~(\ref{d.C_isos})) begin
gives the desired exactness of the top row of~(\ref{d.sharp}).
\endproof

The diagonal chains of donors and receptors which we followed
in the above proof fulfill the promise that ``long''
connections would be reduced to composites of ``short'' ones.
The proof of the next lemma continues this theme.

\begin{lemma}[{Snake Lemma, \cite[p.\,23]{HB}, \cite{2-sq}, \cite[p.\,158]{SL.Alg}, \cite[p.\,50]{McL}}]\label{L.snake}
If, in the commuting
diagram at left below, both rows are exact, and we append
a row of kernels and a row of cokernels to the vertical
maps, as in the diagram at right,
\begin{equation}\begin{minipage}[c]{30pc}\label{d.snake}
\begin{picture}(192,110)

\multiput(0,35)(32,0){2}{
	\multiput(39,0)(0,32){2}{
		\put(0,0){\vector(1,0){18}}}}
\multiput(32,60)(32,0){3}{\vector(0,-1){18}}
\put(-3,32){$0$}\put(7,35){\vector(1,0){18}}
\put(103,67){\vector(1,0){18}}\put(125,64){$0$}
\put(25,32){$Y_1$}\put(57,32){$Y_2$}\put(89,32){$Y_3$}
\put(23,64){$X_1$}\put(55,64){$X_2$}\put(87,64){$X_3$}

\multiput(200,3)(32,0){2}{
	\multiput(39,0)(0,32){4}{
		\put(0,0){\vector(1,0){17}}}}
\multiput(198,28)(0,32){3}{
\multiput(32,0)(32,0){3}{\vector(0,-1){18}}
}
\put(197,32){$0$}\put(207,35){\vector(1,0){18}}
\put(303,67){\vector(1,0){18}}\put(325,64){$0$}
\put(224,00){$C_1$}\put(256,00){$C_2$}\put(288,00){$C_3$}
\put(225,32){$Y_1$}\put(257,32){$Y_2$}\put(289,32){$Y_3$}
\put(223,64){$X_1$}\put(255,64){$X_2$}\put(287,64){$X_3$}
\put(223,95){$K_1$}\put(255,95){$K_2$}\put(287,95){$K_3$}
\put(222,9){\put(0,0){\bx}
\multiput(0,0)(16,16){5}{\dottedline{2}(6,6)(14,14)\put(16,16){\bx}}}
\end{picture}
\end{minipage}\end{equation}
then those two rows fit together into an exact sequence
\begin{equation}\begin{minipage}[c]{30pc}\label{d.snake_seq}
\begin{picture}(400,25)
\multiput(040,13)(50,0){5}{\vector(1,0){25}}
\put(020,10){$K_1$}
\put(070,10){$K_2$}
\put(120,10){$K_3$}
\put(169,10){$C_1$}
\put(219,10){$C_2$}
\put(269,10){$C_3.$}
\end{picture}
\end{minipage}\end{equation}
\end{lemma}

\proof
We extend the right diagram of~(\ref{d.snake}) to a double complex
by attaching a kernel $X_0$ to the second row and a cokernel $Y_4$
to the third, and filling in zeroes everywhere else.
In this complex, the three columns shown
in~(\ref{d.snake}) are exact, and we have
horizontal exactness at $X_1,$ $X_2,$ $X_3,$ $Y_1,$ $Y_2$ and $Y_3.$

The exactness of~(\ref{d.snake_seq}) at $K_2,$ i.e., the triviality
of $K_2\hor,$ now follows from the following isomorphisms
(the first a case of Corollary~\ref{C.inter-iso}, the
next four of Corollary~\ref{C.extra-iso};
cf.~second equation of~(\ref{d.AB_isos})):
\begin{equation}\begin{minipage}[c]{27pc}\label{d.K2_exact}
$K_2\hor\ \cong\ {K_2}_\bx\ \cong\ ^\bx X_2\ \cong\ {X_1}_\bx
\ \cong\ ^\bx\,Y_1\ \cong\ 0_\bx\ =\ 0.$
\end{minipage}\end{equation}
Exactness at $C_2$ is shown similarly.

We now want to find a connecting map $K_3\to C_1$
making~(\ref{d.snake_seq}) exact at these two objects.
This is equivalent to an isomorphism
between $\r{Cok}(K_2\to K_3)=K_3\hor$ and
$\r{Ker}(C_1\to C_2)=C_1\hor.$
And indeed, such an isomorphism is given by the composite
\begin{equation}\begin{minipage}[c]{30pc}
$K_3\hor \cong\ {K_3}_\bx\ \cong\ {^\bx}X_3
\ \cong\ {X_2}_\bx\ \cong\ {^\bx}Y_2\ \cong\ {Y_1}_\bx\ \cong
\ {^\bx}C_1\ \cong\ C_1\hor$
\end{minipage}\end{equation}
of two intramural and five extramural maps
shown as the dotted path in~(\ref{d.snake}), which are again
isomorphisms by Corollaries~\ref{C.inter-iso} and~\ref{C.extra-iso}.
\endproof

The next result, whose proof by the same method we leave to
the reader, establishes isomorphisms between infinitely many pairs of
homology objects in a double complex bordered by zeroes,
either in two parallel, or two perpendicular directions.
Before stating it, we need to make some choices about indexing.

\begin{convention}\label{Cn.cochain}
When the objects of a double complex are indexed by numerical
subscripts, the first subscript will specify the row and the
second the column, and these will increase downwards,
respectively, to the right \textup{(}as in the numbering of the
entries of a matrix; but \emph{not} as in the standard coordinatization
of the $\!(x,y)\!$-plane\textup{)}.
\end{convention}

Since our arrows also point downward and to the right,
our complexes will be double \emph{cochain} complexes;
i.e., the boundary
morphisms will go from lower- to higher-indexed objects.
However, we will continue to call
the constructed objects ``homology objects'', rather
than ``cohomology objects''.

Here is the promised result, which the reader can easily prove
by the method used for Lemmas~\ref{L.sharp} and~\ref{L.snake}.

\begin{lemma}\label{L.2borders}
If in the left-hand complex below, all rows but the first
row shown \textup{(}the row of $\!A_{0,r}\!$'s\textup{)}
and all columns but the first column shown
\textup{(}the column of $\!A_{r,0}\!$'s\textup{)},
are exact, then the homologies of
the first row and the first column are isomorphic:
$A_{0,r}\hor\cong A_{r,0}\ver.$
\textup{(}And analogously for a complex bordered by
zeroes on the \emph{bottom} and the \emph{right}.\textup{)}

If in the right-hand commuting
diagram below, all columns are exact, and all rows
but the first and last are exact, then the homologies of those two rows
agree with a shift of $n-m-1:\ A_{m,r}\hor\cong A_{n,r-n+m+1}\hor.$
\textup{(}And analogously for a complex bordered on
the left and right by zeroes.\textup{)}
\begin{equation}\begin{minipage}[c]{35pc}\label{d.2borders}

\begin{picture}(420,190)
\multiput(20,35)(0,36){4}{
	\put(-3,-3){$0$}\multiput(10,0)(36,0){5}{\vector(1,0){13}}}
\multiput(52,182)(36,0){4}{\put(-4,-6){$0$}
	\multiput(0,-13)(0,-36){5}{\vector(0,-1){16}}}
\put(43,32){
\put(000,000){$A_{30}$}
\put(000,036){$A_{20}$}
\put(000,072){$A_{10}$}
\put(000,108){$A_{00}$}

\put(036,000){$A_{31}$}
\put(036,036){$A_{21}$}
\put(036,072){$A_{11}$}
\put(036,108){$A_{01}$}

\put(072,000){$A_{32}$}
\put(072,036){$A_{22}$}
\put(072,072){$A_{12}$}
\put(072,108){$A_{02}$}

\put(108,000){$A_{33}$}
\put(108,036){$A_{23}$}
\put(108,072){$A_{13}$}
\put(108,108){$A_{03}$}
}

\multiput(230,44)(0,105){2}{
\put(0,0){\vector(1,0){11}}
\put(41,0){\vector(1,0){7}}
\put(75,0){\vector(1,0){10}}
\put(110,0){\vector(1,0){10}}
\put(146,0){\vector(1,0){14}}
}
\multiput(227,76)(0,38){2}{\multiput(0,0)(36,0){5}{\vector(1,0){16}}}
\multiput(254,0)(36,0){4}{
\put(-4,08){$0$}
\put(0,034){\vector(0,-1){16}}
\put(0,069){\vector(0,-1){16}}
\put(-3.5,88){$\vdots$}

\put(0,139){\vector(0,-1){16}}
\put(0,172){\vector(0,-1){16}}
\put(-4,177){$0$}
}
\put(240,145){$A_{m,-1}$}
\put(279,145){$A_{m,0}$}
\put(315,145){$A_{m,1}$}
\put(351,145){$A_{m,2}$}
\multiput(254,76)(0,38){2}{\multiput(0,0)(36,0){4}{\middot}}
\put(240,40){$A_{n,-1}$}
\put(279,40){$A_{n,0}$}
\put(315,40){$A_{n,1}$}
\put(351,40){$A_{n,2}$}

\end{picture}
\end{minipage}\end{equation}
\endproof
\end{lemma}

In contrast to the first of these results, if we form
a double complex bordered on the \emph{top} and the \emph{right}
(or on the \emph{bottom} and the \emph{left}) by zeroes,
and we again assume all rows and columns exact except those adjacent
to the indicated row and column of zeroes, there
will in general be no relation between their homologies.
For a counterexample, one can take a double complex in which all
objects are zero except for a
``staircase'' of isomorphic objects
running upward to the right till it
hits one of the ``borders''.
The place where it hits that border will be the only place
where a nonzero homology object occurs.
(If one tries to construct a similar counterexample to the first
assertion of the above lemma by running a staircase of isomorphisms
upward to the
\emph{left}, one finds that the resulting array of objects
and morphisms is not a commutative diagram.)

\section{Weakly bounded double complexes}\label{S.weakbd}

Before exploring uses of the full statement of the
Salamander Lemma, it will be instructive to consider a mild
generalization of our last result.
Suppose that as in the right-hand diagram of~(\ref{d.2borders})
we have a double complex with exact columns, bounded above the
$\!m\!$-th row and below the $\!n\!$-th row by zeros.
But rather than assuming exactness in all but the $\!m\!$-th and
$\!n\!$-th rows,
let us assume it in all rows but the $\!i\!$-th and $\!j\!$-th,
for some $i$ and $j$ with $m\leq i<j\leq n.$
I claim it will still be true that the homologies of these rows
agree up to a shift:
\begin{equation}\begin{minipage}[c]{30pc}\label{d.ij-hor}
$A_{i,r}\hor\ \cong\ A_{j,r-j+i+1}\hor.$
\end{minipage}\end{equation}

Indeed, first note that by composing extramural isomorphisms
as in the preceding section, we get
\begin{equation}\begin{minipage}[c]{30pc}\label{d.ij-bx}
${A_{i,r}}_\bx\ \cong\ ^\bx A_{j,r-j+i+1}.$
\end{minipage}\end{equation}
So the problem is to strengthen Corollary~\ref{C.inter-iso}
to show that the objects of~(\ref{d.ij-bx})
are isomorphic (by the intramural maps) to their
counterparts in~(\ref{d.ij-hor}).
The required generalization of Corollary~\ref{C.inter-iso} is
quite simple.

\begin{corollary}\label{C.inter-iso+}
Suppose $A$ is an object of a double complex, and the nearby donor and
receptor objects marked ``\put(4,4){\circle{4}}\kern6pt''
in one of the diagrams below are zero.
\begin{equation}\begin{minipage}[c]{30pc}\label{d.iso+}
\begin{picture}(420,105)
\multiput(40,69)(200,0){2}{
\put(0,0){\middot}
\multiput(7,0)(32,0){2}{\vector(1,0){18}}
\put(32,32){\middot}
\multiput(32,25)(0,-34){2}{\vector(0,-1){18}}
\put(26,-3){$A$}
\put(32,-32){\middot}
\put(64,0){\middot}
}
\multiput(0,5)(168,-32){2}{
\multiput(80,88)(16,-16){2}{\circle{4}}
}
\put(47,21){$^\bx A\cong A\ver$}
\put(45,5){$A\hor\cong A_\bx$}
\put(248,21){$^\bx A\cong A\hor$}
\put(250,5){$A\ver\cong A_\bx$}
\end{picture}
\end{minipage}\end{equation}

Then the two intramural isomorphisms indicated below that diagram hold.
\end{corollary}

\proof
To get the first isomorphism of the first
diagram, apply the Salamander Lemma to
the arrow coming vertically into $A;$ to get the second,
apply it to the arrow coming horizontally out of $A.$
In the second diagram, similarly apply it to the two
arrows bearing the ``\put(4,4){\circle{4}}\kern8pt''s.
\endproof

Now in the situation of the first paragraph of this
section, our double complex is
exact horizontally above the $\!i\!$-th row, and vertically
everywhere, so we
can use Corollary~\ref{C.extra-iso} to connect any donor
above the $\!i\!$-th row, or any receptor at or above
that row, to a donor or receptor above the
$\!m\!$-th row, proving it zero.
Corollary~\ref{C.inter-iso+} then shows
the left-hand side of~(\ref{d.ij-hor})
to be isomorphic to the left-hand side of~(\ref{d.ij-bx}).
Similarly, using exactness below the $\!j\!$-th row,
and vanishing below the $\!n\!$-th, we find that
the right-hand sides of those displays are isomorphic.
Thus,~(\ref{d.ij-bx}) yields~(\ref{d.ij-hor}), as desired.

What if we have a double complex in which all columns, and all but
the $\!i\!$-th and $\!j\!$-th rows, are exact, but we do not assume
that all but finitely many rows are zero?
Starting from the receptor at any object of the $\!i\!$-th row,
we can still get an infinite chain of isomorphisms going upward
and to the right:
\begin{equation}\begin{minipage}[c]{30pc}\label{d.stairs}
$^\bx A_{i,s}\ \cong\ {A_{i-1,s}}_\bx\ \cong\ ^\bx A_{i-1,s+1}\ \cong
\ {A_{i-2,s+1}}_\bx\ \cong\ \dots\,,$
\end{minipage}\end{equation}
but we can no longer assert that the common value is zero; and
similarly below the $\!j\!$-th row.
However, there are certainly weaker hypotheses than the one we were
using above that \emph{will} allow us to say this common value
is $0;$ e.g., the existence of zero quadrants (rather than
half-planes) on the upper right and lower left.
Let us make, still more generally,

\begin{definition}\label{D.weaklybdd}
A double complex $(A_{r,s})$ will be called \emph{weakly bounded}
if for every $r$ and $s,$ there exists a \emph{positive} integer $n$
such that $^\bx A_{r-n,s+n}$ or ${A_{r-n-1,s+n}}_\bx$ is zero,
and also a \emph{negative} integer $n$ with the same property.
\end{definition}

The above discussion now yields the first statement of the next
corollary; the final statement is seen to hold by a similar argument.

\begin{corollary}[to proof of Lemma~\ref{L.2borders}]\label{C.i+j}
Let $(A_{r,s})$ be a \emph{weakly bounded} double complex.

If all columns are exact, and all rows but the $\!i\!$-th and
$\!j\!$-th are exact, where $i<j,$ then the
homologies of these
rows are isomorphic with a shift: $A_{i,r}\hor\cong A_{j,r-j+i+1}\hor.$
The analogous statement holds if all rows and all but two columns are
exact.

If all rows but the $\!i\!$-th, and all columns but the $\!j\!$-th
are exact $(i$ and $j$ arbitrary\textup{)}, then the $\!i\!$-th
row and $\!j\!$-th column have isomorphic homologies:
$A_{i,r}\hor\cong A_{r-j+i,j}\ver.$\endproof
\end{corollary}

To see that the above corollary fails without the hypothesis
of weak boundedness, consider again a double complex that is
zero except for a ``staircase'' of copies of a nonzero object
and identity maps between them:
\begin{equation}\begin{minipage}[c]{35pc}\label{d.staircase}
\begin{picture}(420,180)

\multiput(100,0)(0,30){5}{
	\multiput(0,0)(30,0){3}{
		\put(7,30){\vector(1,0){16}}
		\put(26,27){$0$}
		\put(30,23){\vector(0,-1){16}}}
	\put(187,30){\vector(1,0){16}}}
\multiput(190,0)(0,30){2}{
	\multiput(0,0)(30,0){3}{
		\put(7,30){\vector(1,0){16}}
		\put(26,27){$0$}
		\put(30,23){\vector(0,-1){16}}}}
\multiput(130,173)(30,0){5}{\vector(0,-1){16}}
\multiput(197,90)(0,30){3}{\vector(1,0){16}}
\multiput(220,83)(30,0){3}{\vector(0,-1){16}}
\multiput(220,90)(30,0){2}{
	\put(07,0){\vector(1,0){16}}
	\put(30,23){\vector(0,-1){16}}
	\put(26,-3){$0$}}
\multiput(220,90)(30,30){2}{
	\multiput(-6,-3)(0,30){2}{$A$}
	\multiput(-1,23)(2,0){2}{\line(0,-1){16}}
	\multiput(7,28)(0,2){2}{\line(1,0){16}}}
\put(280,150){\put(-6,-3){$A$}
	\multiput(-1,23)(2,0){2}{\line(0,-1){16}}}
\multiput(216,147)(60,-30){2}{$0$}
\multiput(220,143)(60,0){2}{\vector(0,-1){16}}
\multiput(227,150)(30,-30){2}{\vector(1,0){16}}

\end{picture}
\end{minipage}\end{equation}
All rows and all columns are then exact except the \emph{row}
containing the lowest ``$A$''.
Considering that row and any other row, we get a
contradiction to the first conclusion of Corollary~\ref{C.i+j}.
Considering that row and any column gives
a contradiction to the final conclusion.

We remark that if, in Corollary~\ref{C.i+j}, we
make the substitution $s=r+i$ in the subscripts, then our
isomorphisms take the forms $A_{i,s-i}\hor\cong A_{j,s-j+1}\hor$ and
$A_{i,s-i}\hor\cong A_{s-j,j}\ver.$
These formulas are more symmetric than those using $r,$
but I find them a little less easy to think about, because the variable
index $s$ never appears alone.
But in later results, Lemmas~\ref{L.ijk} and~\ref{L.ijk_alt},
where the analog
of the $\!r\!$-indexing would be messier than it is here,
we shall use the analog of this $\!s\!$-indexing.

\section{Long exact sequences}\label{S.long}

At this point it would be easy to apply Lemma~\ref{L.salamander}
and Corollaries~\ref{C.extra-iso} and~\ref{C.inter-iso}
to give a quick construction of the long exact sequence of
homologies associated with a short exact sequences of complexes;
the reader may wish to do so for him or her self.
But we shall find it more instructive to examine how the
six-term exact ``salamander'' sequences we have associated with the
arrows of a double complex
link together under various weaker hypotheses, and see that the
above long exact sequence is the simplest interesting
case of some more general phenomena.

Let $B$ be any object of a double complex, with some neighboring objects
labeled as follows.
\begin{equation}\begin{minipage}[c]{12pc}\label{d.fig8}
\begin{picture}(120,120)
\put(04,107){$D$} \put(54,107){$E$}
\put(04,57){$A$} \put(54,57){$B$} \put(104,57){$C$}
\put(54,07){$F$} \put(104,07){$G\,,$}
\multiput(0,0)(-50,50){2}{
	\multiput(60,50)(50,0){2}{\vector(0,-1){28}}
	\multiput(71,10)(0,50){2}{\vector(1,0){28}}
	}
\end{picture}\end{minipage}\end{equation}
and let us consider the six-term exact sequences associated
with the four arrows into and out of $B.$
These piece together as in the following diagram, where the central
square and each of the four triangular wedges commute:
\begin{equation}\begin{minipage}[c]{32pc}\label{d.4_6s}
\begin{picture}(370,370)

\put(154,357){$D_\bx$}
\put(154,307){$E\ver$}
\put(154,257){$E_\bx$}
\put(150,207){$^\bx B$}
\put(154,157){$B\ver$}

\put(202,207){$B\hor$}
\put(202,157){$B_\bx$}
\put(200,107){$^\bx C$}
\put(204,057){$C\hor$}
\put(200,007){$^\bx G$}

\multiput(160,199)(50,-150){2}{
	\multiput(0,0)(0,50){4}{\vector(0,-1){28}}}

\put(000,207){$D_\bx$}
\put(048,207){$A\hor$}
\put(101,207){$A_\bx$}

\put(250,157){$^\bx F$}
\put(300,157){$F\ver$}
\put(350,157){$^\bx G$}

\multiput(020,210)(150,-50){2}{
	\multiput(0,0)(50,0){4}{\vector(1,0){28}}}

\qbezier(172,260)(210,260)(210,222)\put(210,219){\vector(0,-1){0}}
\qbezier(225,210)(260,210)(260,172)\put(260,169){\vector(0,-1){0}}

\qbezier(107,200)(107,161)(150,161)\put(150,161){\vector(1,0){0}}
\qbezier(160,150)(160,111)(200,111)\put(200,111){\vector(1,0){0}}

\end{picture}\end{minipage}\end{equation}

We now note what happens if~(\ref{d.fig8})
is vertically or horizontally exact at $B.$

\begin{lemma}\label{L.Bver=0}
Suppose in~\textup{(\ref{d.fig8})} that $B\ver=0,$ or that $B\hor=0,$
or more generally, that the intramural map $^\bx B\to B_\bx$ is zero.
Then the following two $\!9\!$-term sequences \textup{(}the
first obtained from the ``left-hand'' and ``bottom'' branches
of~\textup{(\ref{d.4_6s})}, the second from the ``top'' and
``right-hand'' branches\textup{)} are exact:
\begin{equation}\begin{minipage}[c]{35pc}\label{d.9-term1}
\begin{picture}(350,27)
\multiput(30,13)(40,0){8}{\vector(1,0){17}}
\put(10,10){$D_\bx$}
\put(50,10){$A\hor$}
\put(90,10){$A_\bx$}
\put(130,10){$^\bx B$}
\put(170,10){$B\hor$}
\put(210,10){$B_\bx$}
\put(250,10){$^\bx C$}
\put(290,10){$C\hor$}
\put(330,10){$^\bx G,$}
\end{picture}
\end{minipage}\end{equation}
\begin{equation}\begin{minipage}[c]{35pc}\label{d.9-term2}
\begin{picture}(350,27)
\multiput(30,13)(40,0){8}{\vector(1,0){17}}
\put(10,10){$D_\bx$}
\put(50,10){$E\ver$\,}
\put(90,10){$E_\bx$}
\put(130,10){$^\bx B$}
\put(170,10){$B\ver$\,}
\put(210,10){$B_\bx$}
\put(250,10){$^\bx F$}
\put(290,10){$F\ver$\,}
\put(330,10){$^\bx G.$}
\end{picture}
\end{minipage}\end{equation}
\end{lemma}
\proof
The exactness of the $\!6\!$-term sequences
of which~(\ref{d.4_6s}) is composed gives the exactness
of~(\ref{d.9-term1}) and~(\ref{d.9-term2}) everywhere
but at the middle terms, $B\hor$ and $B\ver.$
That the composite map through that middle term equals
zero is, in each case, our hypothesis on the intramural
map $^\bx B\to B_\bx.$
That, conversely, the kernel of the map out of that middle object is
contained in the image of the map going into it can be
seen from~(\ref{d.4_6s}); e.g., in the case of~(\ref{d.9-term1}),
we see from~(\ref{d.4_6s}) that the kernel of the
map $B\hor\to B_\bx$ is the image of the map $E_\bx\to B\hor,$
and that map factors through the map $^\bx B\to B\hor.$
(This is the commutativity of the topmost
triangular wedge in~(\ref{d.4_6s}).
A dual proof can be gotten using the right-hand wedge.)
\endproof

In noting applications of the above result, we shall,
for brevity, restrict ourselves to~(\ref{d.9-term1}); the
corresponding consequences of~(\ref{d.9-term2}) follow by symmetry.
Of the alternative hypotheses of Lemma~\ref{L.Bver=0}, the condition
$B\hor=0$ makes~(\ref{d.9-term1}) degenerate, while the more
general statement that $^\bx B\to B_\bx$ is zero does not
correspond to any condition
in the standard language of homological algebra; so in the
following corollary, we focus mainly on the condition $B\ver=0.$

\begin{corollary}\label{C.1-row-long}
If in a double complex, a piece of which is labeled as
in~\textup{~(\ref{d.fig8})}, the \emph{vertical}
homologies are zero for all objects in the
{\em row} $\dots\to A\to B\to C\to\dots$
\textup{(}or more generally, if the intramural map
from receptor to donor is zero for each object of
that row\textup{)}, then the following sequence
of horizontal homology objects, donors
and receptors, and intramural and extramural maps,
is exact:
\begin{equation}\begin{minipage}[c]{35pc}\label{d.1-row-long}
\begin{picture}(400,27)
\multiput(28,13)(35,0){10}{\vector(1,0){15}}
\put(8,10){$\cdots$}
\put(44,10){$^\bx A$}
\put(78,10){$A\hor$}
\put(115,10){$A_\bx$}
\put(149,10){$^\bx B$}
\put(182,10){$B\hor$}
\put(220,10){$B_\bx$}
\put(255,10){$^\bx C$}
\put(288,10){$C\kern-.1em\hor$}
\put(325,10){$C_\bx$}
\put(360,10){$\cdots\,.$}
\end{picture}
\end{minipage}\end{equation}
\end{corollary}

\proof
At each object of the indicated row of~(\ref{d.fig8}),
write down the exact sequence corresponding
to~(\ref{d.9-term1}), leaving off
the first and last terms.
The remaining parts of
these sequences overlap, giving~(\ref{d.1-row-long}).
\endproof

When the vertical homology in our double complex is \emph{everywhere}
zero,
the exact sequences~(\ref{d.1-row-long}) arising from
successive rows are linked, at every third position, by
isomorphisms given by Corollary~\ref{C.extra-iso}, as described in

\begin{corollary}\label{C.linked_long}
If in a double complex
\begin{equation}\begin{minipage}[c]{21pc}\label{d.A-Z}
\begin{picture}(210,150)

\put(024,27){$W$}
\put(054,27){$X$}
\put(084,27){$Y$}
\put(114,27){$Z$}

\put(024,57){$P$}
\put(054,57){$Q$}
\put(084,57){$R$}
\put(114,57){$S$}

\put(054,87){$K$}
\put(084,87){$L$}
\put(114,87){$M$}
\put(144,87){$N$}

\put(084,117){$A$}
\put(114,117){$B$}
\put(144,117){$C$}
\put(174,117){$D$}

\multiput(0,30)(30,30){4}{
\multiput(0,0)(30,0){3}{
\put(7,0){\vector(1,0){16}}\put(30,23){\vector(0,-1){16}}}
\put(97,0){\vector(1,0){16}}
}
\multiput(0,60)(30,30){3}{
\put(7,0){\vector(1,0){16}}\put(30,23){\vector(0,-1){16}}
\put(120,-07){\vector(0,-1){16}}
}
\multiput(0,0)(30,0){4}{
\put(30,23){\vector(0,-1){16}}
}
\put(127,30){\vector(1,0){16}}
\end{picture}
\end{minipage}\end{equation}
all \emph{columns} are exact, then the \emph{rows} induce long
exact sequences, which are linked by isomorphisms:
\begin{equation}\begin{minipage}[c]{33pc}\label{d.linked_long}
\begin{picture}(385,150)
\put(024,27){$^\bx W$}
\put(058,27){$W\kern-.25em\hor$}
\put(095,27){$W_\bx$}
\put(128,27){$^\bx X$}
\put(165,27){$X\kern-.20em\hor$}
\put(199,27){$X_\bx$}
\put(234,27){$^\bx Y$}
\put(270,27){$Y\kern-.25em\hor$}
\put(307,27){$Y_\bx$}
\put(339,27){$^\bx Z$}

\put(026,57){$P_\bx$}
\put(057,57){$^\bx Q$}
\put(093,57){$Q\hor$}
\put(131,57){$Q_\bx$}
\put(163,57){$^\bx R$}
\put(199,57){$R\hor$}
\put(235,57){$R_\bx$}
\put(268,57){$^\bx S$}
\put(304,57){$S\hor$}
\put(341,57){$S_\bx$}

\put(024,87){$K\kern-.20em\hor$}
\put(059,87){$K_\bx$}
\put(091,87){$^\bx L$}
\put(128,87){$L\kern.05em\hor$}
\put(163,87){$L_\bx$}
\put(195,87){$^\bx M$}
\put(231,87){$M\kern-.15em\hor$}
\put(268,87){$M_\bx$}
\put(303,87){$^\bx N$}
\put(338,87){$N\kern-.15em\hor$}

\put(023,117){$^\bx A$}
\put(058,117){$A\hor$}
\put(093,117){$A_\bx$}
\put(127,117){$^\bx B$}
\put(162,117){$B\hor$}
\put(198,117){$B_\bx$}
\put(231,117){$^\bx C$}
\put(268,117){$C\hor$}
\put(305,117){$C_\bx$}
\put(338,117){$^\bx D$}

\multiput(0,30)(0,30){4}{
\multiput(0,0)(35,0){11}{
\put(8,0){\vector(1,0){14}}}}

\multiput(105,22)(0,90){2}{
\multiput(0,0)(105,0){3}{
\multiput(-3,0)(2,0){2}{\line(0,-1){14}}
}}

\multiput(35,52)(0,90){2}{
\multiput(0,0)(105,0){4}{
\multiput(-3,0)(2,0){2}{\line(0,-1){14}}
}}

\multiput(70,82)(105,0){3}{
\multiput(-1,0)(2,0){2}{\line(0,-1){14}}
}
\end{picture}
\endproof
\end{minipage}\end{equation}
\end{corollary}

Note that in these long exact sequences, the classical homology objects
form every third term -- the terms of~(\ref{d.linked_long})
that are not connected either above or below by isomorphisms.

Suppose now that in~(\ref{d.A-Z}), in addition to all columns
being exact, some row is exact.
This means that in the system of long exact
sequences~(\ref{d.linked_long}),
the corresponding row will have every third term zero; and so
the maps connecting the remaining terms will be isomorphisms:
\begin{equation}\begin{minipage}[c]{20pc}\label{d.->&=}
\begin{picture}(240,100)

\multiput(0,20)(0,20){2}{
\multiput(0,0)(20,0){10}{
\put(4,0){\vector(1,0){10}}\put(21,0){\middot}}
\put(204,0){\vector(1,0){10}}
}
\put(0,80){
\multiput(0,0)(20,0){10}{
\put(4,0){\vector(1,0){10}}\put(21,0){\middot}}
\put(204,0){\vector(1,0){10}}
}

\multiput(-1,15)(60,0){3}{
\multiput(20,20)(20,20){4}{
\multiput(0,0)(2,0){2}{\line(0,-1){10}}
}
\multiput(60,0)(2,0){2}{\line(0,-1){10}}
}

\multiput(199,35)(2,0){2}{\line(0,-1){10}}

\multiput(19,95)(2,0){2}{\line(0,-1){10}}

\multiput(40,60)(60,0){3}{
\put(1,0){\middot}
\put(5,0){\multiput(0,0)(0,2){2}{\line(1,0){10}}}
\put(21,0){\middot}
}
\end{picture}
\end{minipage}\end{equation}

We see that these, together with the vertical isomorphisms,
tie together the \emph{preceding} and
\emph{following} exact sequence to give a
system essentially like~(\ref{d.linked_long}), except for a
horizontal shift by one step.
If $n$ successive rows of the double complex are exact, we get a
similar diagram with a shift by $n$ steps.

If all rows are exact above a certain point, then we get infinite
chains of isomorphisms going upward and to the right.
If the complex is also weakly bounded (Definition~\ref{D.weaklybdd}),
the common value along those chains will be zero; hence
every third term of the long exact sequence corresponding
to the top \emph{nonexact} row of our double complex
will be zero, so we again have
isomorphisms between pairs of remaining terms; though not the
same pairs as before: in each isomorphic pair,
one of the members is now a horizontal homology object.
If we regard classical
homology objects as more interesting than donors and
receptors, we may use these isomorphisms and the isomorphisms
joining this row to the next to insert these homology
objects in that row, in place of all the receptors.

For instance, if all rows of~(\ref{d.A-Z}) above the top one
shown are exact, and the complex is weakly bounded above,
then in the top two rows of~(\ref{d.linked_long}) we get
\raisebox{-2.8em}{\kern-.05em
\begin{picture}(140,68)
\multiput(4,30)(30,0){5}{\vector(1,0){12}}
\put(15,27){$K\kern-.2em\hor$}
\put(45,27){$K_\bx$}
\put(75,27){$^\bx L$}
\put(105,27){$L\hor$}

\put(45,57){$A\hor$}
\put(75,57){$A_\bx$}

\multiput(55,20)(2,0){2}{\line(0,-1){10}}
\multiput(85,50)(2,0){2}{\line(0,-1){10}}
\multiput(65,59)(0,2){2}{\line(1,0){10}}
\end{picture}},
which we can rewrite
\raisebox{-2.0em}{
\begin{picture}(140,25)
\multiput(4,30)(30,0){5}{\vector(1,0){12}}
\put(15,27){$K\kern-.2em\hor$}
\put(45,27){$K_\bx$}
\put(75,27){$A\hor$}
\put(105,27){$L\hor$}
\multiput(55,20)(2,0){2}{\line(0,-1){10}}
\end{picture}}.

Of course, if, say, the second row of~(\ref{d.A-Z})
(unlike the first) happens to
be exact, then the objects $K\hor,$ $L\hor$ etc.\ in the
above exact sequence are zero, allowing us to
pull the homology objects from the first
row down yet another step, and insert them into the long exact
sequence arising from the third row.
Continuing in this way as long as we find exact rows in~(\ref{d.A-Z}),
we get a linked system of long exact sequences, of which
the top sequence has, as \emph{two} out of every three terms,
horizontal homology objects, and arises from the two highest
non-exact rows of~(\ref{d.A-Z}) (assuming there are at least two).
The obvious analogous situation holds
if, instead, all rows \emph{below} some point are exact.

If our original double complex
has only three nonexact rows, then we can see that, working
in this way from both ends, we get a single long
exact sequence with horizontal homology objects for all its terms:

\begin{lemma}\label{L.ijk}
Suppose we are given a weakly bounded double complex,
with objects $A_{h,r},$ all columns exact, and all rows exact
except the $\!i\!$-th, $\!j\!$-th and $\!k\!$-th, where
$i<j<k.$
Then we get a long exact sequence
\begin{equation}\begin{minipage}[c]{35pc}\label{d.ijk}
$\cdots{\to}
A_{i,\,s-i-1}\hor\,{\to}\,
A_{j,\,s-j}\hor\,{\to}\,
A_{k,\,s-k+1}\hor\,{\to}\,
A_{i,\,s-i}\hor\,{\to}\,
A_{j,\,s-j+1}\hor\,{\to}\,
A_{k,\,s-k+2}\hor{\to}\,{\cdots}.$
\end{minipage}\end{equation}
\textup{(}Regarding the indexing, cf.~last paragraph of
\S\ref{S.weakbd}.\textup{)}
\endproof
\end{lemma}

What if we have four rather than three non-exact rows (again in
a weakly bounded double complex with exact columns)?
Assuming for concreteness that our double complex
is~(\ref{d.A-Z}), and that all rows but
the four shown there are exact, we find
that~(\ref{d.linked_long}) collapses to
\begin{equation}\begin{minipage}[c]{32pc}\label{d.2_lined}
\begin{picture}(385,51)

\put(023,07){$W\kern-.25em\hor$}
\put(061,07){$^\bx Q$}
\put(095,07){$Q\hor$}
\put(130,07){$X\kern-.25em\hor$}
\put(166,07){$^\bx R$}
\put(200,07){$R\kern-.1em\hor$}
\put(235,07){$Y\hor$}
\put(271,07){$^\bx S$}
\put(305,07){$S\hor$}
\put(340,07){$Z\hor$}

\put(023,42){$K\kern-.25em\hor$}
\put(061,42){$K_\bx$}
\put(095,42){$A\hor$}
\put(130,42){$L\hor$}
\put(166,42){$L_\bx$}
\put(200,42){$B\hor$}
\put(234,42){$M\kern-.25em\hor$}
\put(270,42){$M_\bx$}
\put(305,42){$C\hor$}
\put(339,42){$N\kern-.25em\hor$}

\multiput(0,11)(0,35){2}{
\multiput(0,0)(35,0){11}{
\put(8,0){\vector(1,0){14}}}}

\multiput(70,35)(105,0){3}{
\multiput(-1,0)(2,0){2}{\line(0,-1){14}}
}

\multiput(-13,07)(390,0){2}{\multiput(0,0)(0,35){2}{$\cdots$}}
\put(400,8){.}

\end{picture}
\end{minipage}\end{equation}

\section{Some rows, and some columns}\label{S.rows+cols}

We have just seen what happens when all columns, and all
but a finite number of rows of a weakly bounded double complex
are exact; the corresponding results hold, of course, when all
rows and all but a finite number of columns are exact.

One can look, more generally, at the situation where
\begin{equation}\begin{minipage}[c]{30pc}\label{d.m+n}
All but $m$ rows, and all but $n$ columns are exact.
\end{minipage}\end{equation}

In this section we shall examine the sort of behavior that this
leads to.
For the first result, I give a formal statement,
Lemma~\ref{L.ijk_alt}, and a sketch of
the proof, of which the reader can check the details,
following the technique of the preceding section.
The proofs of the remaining results discussed use the same ideas.

The case $m+n\leq 2$ of~(\ref{d.m+n})
is covered by Corollary~\ref{C.i+j} above.
The case $m+n=3$ is covered (up to row-column reversal) by
Lemma~\ref{L.ijk}, together with

\begin{lemma}\label{L.ijk_alt}
Suppose we have a weakly bounded double complex with objects $A_{h,r},$
all rows exact but the $\!i\!$-th and $\!j\!$-th, where $i<j,$ and
all columns exact but the $\!k\!$-th.
Then we have a long exact sequence
\begin{equation}\begin{minipage}[c]{35pc}\label{d.ijk_alt}
$\cdots\,\to
\ A_{i,s-i-1}\hor\,{\to}
\ A_{j,s-j}\hor\,{\to}
\ A_{s-k,k}\,\ver\,{\to}
\ A_{i,s-i}\hor\,{\to}
\ A_{j,s-j+1}\hor\,{\to}
\ A_{s-k+1,k}\,\ver\,{\to}\,\cdots\ .$
\end{minipage}\end{equation}
\end{lemma}

\subsubsection*{\textsc{Sketch of Proof}}
Write out the salamander exact sequences corresponding
to the \emph{horizontal} arrow out of $A_{i,r}$ for each $r<k,$
to the \emph{vertical} arrow out of $A_{h,k}$ for $h=i,\dots,j-1,$
and to the \emph{horizontal} arrow out of $A_{j,r}$ for each $r\geq k.$

Except where we come to a corner, these exact sequences piece together
(due to exactness of all other rows and columns)
as in~(\ref{d.1-row-long}).
When we do turn a corner, we get a different sort of piecing together;
e.g., if we take the maps $A_{i,k-1}\to A_{i,k}\to A_{i+1,k}$
for the arrows $A\to B\to F$ of~(\ref{d.fig8}), then the
$E_\bx$ and $^\bx C$ of~(\ref{d.fig8}) are both zero,
due to weak boundedness, so
that in~(\ref{d.4_6s}), the two horizontal exact sequences
collapse into one $\!8\!$-term exact sequence.
So the path of arrows in our double complex
described in the preceding paragraph leads to a single
long exact sequence of homology objects, donors, and receptors.

For each donor or receptor in this
sequence, we now use a string of extramural isomorphisms
(consequences of Corollary~\ref{C.extra-iso} and
the exactness of all but our three exceptional rows and columns)
to connect it with a receptor or
donor at an object of one of the other two non-exact rows or columns.
(In each case, there is only one direction we can go by extramural
isomorphisms from our donor or receptor object, without crossing
the non-exact row or column we are on, and this
indeed leads to an object of another non-exact row or column.)
We know from weak boundedness that the donor and receptor
objects on the \emph{other} side of the row or column we have
arrived at are zero, and so conclude by Corollary~\ref{C.inter-iso+}
that the receptor or donor we have reached is isomorphic to
a vertical or horizontal homology object in that row or column.

Thus, we get an exact sequence in which all objects are
vertical or horizontal homology objects.
\endproof

So far, the general case of~(\ref{d.m+n})
has given results as nice as when $m$ or $n$ is $0.$
But now consider $m+n=4.$
We saw in~(\ref{d.2_lined}) what happens
when $n=0;$ let us compare this with the case
of a double complex with three not necessarily exact
rows and one not necessarily exact column, such as the following
(where we have darkened the arrows in the not necessarily
exact rows and columns).
\begin{equation}\begin{minipage}[c]{25pc}\label{d.3+1}
\begin{picture}(360,190)

\multiput(100,30)(0,30){6}{
	\multiput(0,0)(30,0){5}{\vector(0,-1){16}}
	\put(60,0){\thicklines\vector(0,-1){16}}
	}
\multiput(101,37)(0,120){2}{
	\multiput(-23,0)(30,0){6}{\vector(1,0){16}}
	\multiput(0,0)(90,0){2}{
	\multiput(0,0)(30,0){2}{\middot}
	}
	}
\multiput(100,67)(0,30){3}{
	\multiput(-23,0)(30,0){6}{\thicklines\vector(1,0){16}}
	}
\put(93,33){
	\put(60,120){$X$}

	\put(0,90){$A$}
	\put(30,90){$B$}
	\put(60,90){$C$}
	\put(90,90){$D$}
	\put(120,90){$E$}

	\put(0,60){$F$}
	\put(30,60){$G$}
	\put(60,60){$H$}
	\put(90,60){$J$}
	\put(120,60){$K$}

	\put(0,30){$L$}
	\put(30,30){$M$}
	\put(60,30){$N$}
	\put(90,30){$P$}
	\put(120,30){$Q$}

	\put(60,0){$Y$}
}
\end{picture}
\end{minipage}\end{equation}

One finds that this double complex leads to a system of four linked
``half-long'' exact sequences.
To the left and to the right, the diagram looks like the
$n=0$ case,~(\ref{d.2_lined}),
but there is a peculiar ``splicing'' in the middle:
\begin{equation}\begin{minipage}[c]{35pc}\label{d.splice}
\begin{picture}(420,110)

\multiput(-8,10)(0,60){2}{\multiput(0,30)(224,-30){2}{
\multiput(0,0)(28,0){8}{\vector(1,0){10}}
}
\multiput(172,0)(56,30){2}{\vector(1,0){34}}
}
\put(-22,33){
\multiput(168,00)(113,0){2}{\multiput(0,0)(0,60){2}{
	\multiput(0,0)(2,0){2}{\line(0,-1){16}}}}
\multiput(30,50)(224,-30){2}{\multiput(0,0)(84,0){3}{
	\multiput(0,0)(2,0){2}{\line(0,-1){26}}}}
\multiput(34,50)(308,-30){2}{\multiput(0,0)(84,0){2}{
	\qbezier(0,0)(2,-6.5)(0,-13)
	\qbezier(0,-13)(-2,-19.5)(0,-26)}}
\multiput(224,30)(2,0){2}{\line(0,-1){16}}
}
\put(196,55){\oval(10,80)[l]}
\put(195,55){\oval(12,84)[l]}
\put(2,7){
\put(0,90){$A_\bx$}
\put(28,90){$X\kern-.1em\ver$}
\put(54,90){$B\hor$}
\put(84,90){$B_\bx$}
\put(112,90){$C\ver$}
\put(139,90){$C\kern-.10em\hor$}
\put(168,90){$^\bx H$}
\put(196,90){$H\ver$}
\put(251,90){$D\kern-.10em\hor$}

\put(140,60){$C\hor$}
\put(194,60){$H\kern-.15em\hor$}
\put(224,60){$H_\bx$}
\put(251,60){$D\kern-.10em\hor$}
\put(280,60){$J\kern-.25em\hor$}
\put(308,60){$J_\bx$}
\put(335,60){$E\kern-.10em\hor$}
\put(364,60){$K\kern-.25em\hor$}
\put(392,60){$K_\bx$}

\put(0,30){$^\bx F$}
\put(28,30){$F\kern-.25em\hor$}
\put(56,30){$L\hor$}
\put(84,30){$^\bx G$}
\put(111,30){$G\kern-.10em\hor$}
\put(138,30){$M\kern-.20em\hor$}
\put(168,30){$^\bx H$}
\put(196,30){$H\kern-.15em\hor$}
\put(252,30){$N\kern-.25em\hor$}

\put(138,0){$M\kern-.20em\hor$}
\put(196,0){$H\ver$}
\put(224,0){$H_\bx$}
\put(251,0){$N\kern-.15em\hor$}
\put(280,0){$N\kern-.1em\ver$}
\put(308,0){$^\bx P$}
\put(335,0){$P\kern-.15em\hor$}
\put(364,0){$Y\kern-.1em\ver$}
\put(392,0){$^\bx Q$}
}
\end{picture}
\end{minipage}\end{equation}

Here is the same diagram, redrawn more smoothly.

\begin{equation}\begin{minipage}[c]{35pc}\label{d.smooth_splice}
\begin{picture}(420,120)
\put(2,6){
\put(-16,50){$^\bx F\cong A_\bx$}
\put(26,90){$X\kern-.1em\ver$}
\put(26,10){$F\kern-.15em\hor$}
\put(58,90){$B\hor$}
\put(58,10){$L\hor$}
\put(68,50){$^\bx G\cong B_\bx$}
\put(110,10){$C\ver$}
\put(110,90){$G\hor$}
\put(142,10){$C\hor$}
\put(142,90){$M\kern-.10em\hor$}
\put(169,50){$^\bx H$}
\put(196,80){$H\ver$}
\put(196,20){$H\kern-.10em\hor$}
\put(222,50){$H_\bx$}
\put(250,90){$D\hor$}
\put(250,10){$N\kern-.10em\hor$}
\put(282,90){$J\kern-.10em\hor$}
\put(282,10){$N\kern-.1em\ver$}
\put(292,50){$^\bx P\cong J_\bx$}
\put(334,10){$E\hor$}
\put(334,90){$P\kern-.05em\hor$}
\put(366,10){$K\kern-.10em\hor$}
\put(366,90){$Y\kern-.1em\ver$}
\put(376,50){$^\bx Q\,{\cong}\,K_\bx$}
}
\put(13,60){
\multiput(0,0)(84,0){2}{\multiput(0,0)(224,0){2}{
\put(1,7){\vector(2,3){17}}
\put(1,-7){\vector(2,-3){17}}
\put(62,33){\vector(2,-3){17}}
\put(62,-33){\vector(2,3){17}}
\put(28,49){\qbezier(0,-2)(10,8)(20,0)\put(23,-3){\vector(1,-1){0}}}
\put(28,-49){\qbezier(0,2)(10,-8)(20,0)\put(23,3){\vector(1,1){0}}}
}}
\put(-23,33){\vector(2,-3){17}}
\put(-23,-33){\vector(2,3){17}}
\put(396,8){\vector(2,3){17}}
\put(396,-8){\vector(2,-3){17}}
\put(176,8){\vector(3,4){10}}
\put(176,-8){\vector(3,-4){10}}
\put(204,23){\vector(3,-4){10}}
\put(204,-23){\vector(3,4){10}}
\put(193,0){
     \qbezier(-49,+48)(-31,+59)(-7,+38)\put(-6,+37){\vector(1,-1){0}}
     \qbezier(-49,-48)(-31,-59)(-7,-38)\put(-6,-37){\vector(1,+1){0}}
     \qbezier(+49,+48)(+31,+59)(+7,+38)\put(52,+45){\vector(1,-1){0}}
     \qbezier(+49,-48)(+31,-59)(+7,-38)\put(52,-45){\vector(1,+1){0}}
}
}
\end{picture}
\end{minipage}\end{equation}
The exact sequences are those chains of arrows which can be
followed without making sharp turns.

We remark that the first step in verifying the exactness of
the sequences in~(\ref{d.splice}), equivalently,~(\ref{d.smooth_splice})
is to check that the following isomorphisms follow from
Corollaries~\ref{C.extra-iso} and~\ref{C.inter-iso+}.\vspace{.3em}
\begin{equation}\begin{minipage}[c]{32pc}\label{d.3+1_isos}
$X\kern-.1em\ver\cong{^\bx X}\cong{^\bx B},$\quad
$C\ver\cong{^\bx C},$\quad
$C\hor\cong C_\bx\,,$\quad
$D\hor\cong D_\bx\cong{^\bx J}\,,$\quad
$^\bx P\cong J_\bx\,,$\\[.5em]
$^\bx Q\cong K_\bx\,,$\ \ \qquad
$E\hor\cong E_\bx\cong{^\bx K},$\qquad\quad
$L\hor\cong{^\bx L}\cong F_\bx\,,$\ \ \qquad
$^\bx F\cong A_\bx\,,$\\[.5em]
$Y\kern-.1em\ver\cong Y_\bx\cong P_\bx,$\quad
$N\kern-.1em\ver\cong N_\bx,$\quad
$N\hor\cong{^\bx N},$\quad
$M\hor\cong{^\bx M}\cong G_\bx,$\quad
$^\bx G\cong B_\bx\,.$\quad
\end{minipage}\end{equation}
Using these, the verification of the exactness conditions is immediate.
(Note that the part of~(\ref{d.smooth_splice}) between
``$^\bx G\cong B_\bx$'' and ``$^\bx P\cong J_\bx$'' is, up to
labeling, just a copy of~(\ref{d.4_6s}), with the top and left nodes
of~(\ref{d.4_6s}) identified, and likewise the bottom and right nodes,
and with substitutions from~(\ref{d.3+1_isos}) made where appropriate.
So all the exactness conditions in that part of the
diagram are immediate.)

The interpolation of some exact rows between the three nonexact
rows of~(\ref{d.3+1}) does not affect the resulting system
of exact sequences~(\ref{d.splice}),~(\ref{d.smooth_splice})
except by a shift of indices.

General values of $m$ and $n$ in~(\ref{d.m+n})
yield systems that, for most
of their length, consist of $m+n-2$ intertwining exact
sequences (cf.~(\ref{d.->&=})), but have finitely many
``splicings''; essentially, one for each object
of the given double complex which lies at the intersection
of a nonexact row and a nonexact column, and does \emph{not} have,
either to its upper right or lower left, a region where
(due to exactness and weak boundedness) all donors
and receptors are zero.
Thus, the one splice in~(\ref{d.splice}) and~(\ref{d.smooth_splice})
comes from the object $H$ of~(\ref{d.3+1}).

If $m=n=2,$ as in
\begin{equation}\begin{minipage}[c]{35pc}\label{d.2+2}
\begin{picture}(360,220)

\put(94,33){
	\put(60,150){$P$}
	\put(90,150){$U$}

	\put(60,120){$Q$}
	\put(90,120){$V$}

	\put(0,90){$A$}
	\put(30,90){$B$}
	\put(60,90){$C$}
	\put(90,90){$D$}
	\put(120,90){$E$}
	\put(150,90){$F$}

	\put(0,60){$H$}
	\put(30,60){$J$}
	\put(60,60){$K$}
	\put(90,60){$L$}
	\put(119,60){$M$}
	\put(150,60){$N$}

	\put(60,30){$S$}
	\put(90,30){$X$}

	\put(61,0){$T$}
	\put(90,0){$Y$}
}

\multiput(100,0)(120,0){2}{\multiput(0,0)(30,0){2}{
	\multiput(0,30)(0,30){7}{\vector(0,-1){16}}
	\multiput(0,37)(0,120){2}{\multiput(0,0)(0,30){2}{\middot}}
	}}
\multiput(160,0)(30,0){2}{
	\multiput(0,30)(0,30){7}{\thicklines\vector(0,-1){16}}
	}
\multiput(100,37)(0,120){2}{\multiput(0,0)(0,30){2}{
	\multiput(-23,0)(30,0){7}{\vector(1,0){16}}
	}}
\multiput(100,97)(0,30){2}{
	\multiput(-23,0)(30,0){7}{\thicklines\vector(1,0){16}}
	}
\end{picture}
\end{minipage}\end{equation}
then the two objects $C$ and $L$ lead to two ``splicings'':
\begin{equation}\begin{minipage}[c]{35pc}\label{d.2-splice}
\begin{picture}(420,130)

\put(2,7){
\put(0,90){$^\bx Q$}
\put(30,90){$Q\ver$}
\put(60,90){$V\ver$}
\put(90,90){$^\bx C$}
\put(119,90){$C\ver$}
\put(180,90){$D\ver$}
\put(208,90){$D\hor$}
\put(270,90){$L\kern-.10em\hor$}
\put(300,90){$L_\bx$}
\put(329,90){$E\kern-.10em\hor$}
\put(358,90){$M\kern-.20em\hor$}
\put(390,90){$M_\bx$}

\put(60,60){$V\ver$}
\put(119,60){$C\kern-.10em\hor$}
\put(150,60){$C_\bx$}
\put(180,60){$D\ver$}
\put(208,60){$D\hor$}
\put(240,60){$^\bx L$}
\put(270,60){$L\ver$}
\put(329,60){$E\kern-.10em\hor$}

\put(0,30){$^\bx B$}
\put(28,30){$B\kern-.05em\hor$}
\put(60,30){$J\kern-.15em\hor$}
\put(90,30){$^\bx C$}
\put(119,30){$C\kern-.10em\hor$}
\put(180,30){$K\kern-.15em\hor$}
\put(210,30){$K\ver$}
\put(270,30){$L\ver$}
\put(300,30){$L_\bx$}
\put(330,30){$S\ver$}
\put(360,30){$X\kern-.1em\ver$}
\put(390,30){$X_\bx$}

\put(60,00){$J\kern-.15em\hor$}
\put(119,00){$C\ver$}
\put(150,00){$C_\bx$}
\put(179,00){$K\kern-.20em\hor$}
\put(210,00){$K\ver$}
\put(240,00){$^\bx L$}
\put(270,00){$L\hor$}
\put(330,00){$S\ver$}
}

\put(120,55){\oval(10,80)[l]\put(-1,0){\oval(12,84)[l]}}
\put(270,55){\oval(10,80)[l]\put(-1,0){\oval(12,84)[l]}}

\multiput(-10,10)(0,60){2}{
	\multiput(0,30)(300,0){2}{
		\multiput(0,0)(30,0){5}{\vector(1,0){12}}}
	\multiput(150,30)(90,0){2}{\vector(1,0){42}}

	\multiput(150,0)(30,0){5}{\vector(1,0){12}}
	\multiput(90,0)(210,0){2}{\vector(1,0){42}}
	\put(16,23){
	\multiput(60,00)(150,0){2}{\multiput(0,0)(120,0){2}{
	\multiput(0,0)(2,0){2}{\line(0,-1){14}}}}
	}
}
\put(8,33){
\multiput(0,50)(390,0){2}{\multiput(0,0)(2,0){2}{\line(0,-1){26}}}
\multiput(90,53)(210,0){2}{\multiput(0,0)(2,0){2}{\line(0,-1){32}}}
\multiput(120,30)(150,0){2}{\multiput(0,0)(2,0){2}{\line(0,-1){14}}}
\multiput(150,23)(90,0){2}{\multiput(0,0)(2,0){2}{\line(0,-1){32}}}
}
\multiput(14,83)(390,0){2}{
	\qbezier(0,0)(2,-6.5)(0,-13)
	\qbezier(0,-13)(-2,-19.5)(0,-26)}
\end{picture}
\end{minipage}\end{equation}

The same diagram in ``smooth'' format (and carried one step
further at each end) is
\begin{equation}\begin{minipage}[c]{35pc}\label{d.smooth_2-splice}
\begin{picture}(420,85)
\put(-10,37){

\put(0,0){$^\bx A\cong\kern-0.1em{^\bx P}$}
\put(35,27){$P\ver$}
\put(35,-27){$A\hor$}
\put(60,27){$U\ver$}
\put(60,-27){$H\kern-.15em\hor$}
\put(60,0){$^\bx Q\cong\kern-0.1em{^\bx B}$}
\put(95,27){$B\hor$}
\put(95,-27){$Q\ver$}
\put(120,27){$J\kern-.15em\hor$}
\put(120,-27){$V\ver$}
\put(135,0){$^\bx C$}
\put(156,20){$C\ver$}
\put(156,-20){$C\kern-.15em\hor$}
\put(175,0){$C_\bx$}
\put(194,27){$D\ver$}
\put(194,-27){$K\kern-.15em\hor$}
\put(215,27){$D\hor$}
\put(215,-27){$K\ver$}
\put(230,0){$^\bx L$}
\put(252,20){$L\hor$}
\put(252,-20){$L\ver$}
\put(270,0){$L_\bx$}
\put(290,27){$S\ver$}
\put(290,-27){$E\hor$}
\put(310,27){$X\kern-.1em\ver$}
\put(310,-27){$M\kern-.15em\hor$}
\put(315,0){$M_\bx\cong X_\bx$}
\put(350,27){$F\hor$}
\put(350,-27){$T\ver$}
\put(375,27){$N\kern-.20em\hor$}
\put(375,-27){$Y\kern-.1em\ver$}
\put(375,0){$N_\bx\cong Y_\bx$}
}

\put(15,40){
\multiput(0,0)(60,0){2}{\multiput(0,0)(255,0){2}{
\put(2,7){\vector(2,3){10}}
\put(2,-7){\vector(2,-3){10}}
\put(44,22){\vector(2,-3){10}}
\put(44,-21){\vector(2,3){10}}
\put(20,33){\qbezier(0,0)(7,8)(15,0)\put(16,-1){\vector(1,-1){0}}}
\put(20,-33){\qbezier(0,0)(7,-8)(15,0)\put(16,1){\vector(1,1){0}}}
}}

\put(-15,22){\vector(2,-3){10}}
\put(-15,-22){\vector(2,3){10}}
\put(375,5){\vector(2,3){10}}
\put(375,-5){\vector(2,-3){10}}
\multiput(125,0)(95,0){2}{
\put(0,5){\vector(2,3){7}}
\put(0,-5){\vector(2,-3){7}}
\put(20,15){\vector(2,-3){7}}
\put(20,-15){\vector(2,3){7}}
}
\multiput(137,0)(95,0){2}{
     \qbezier(-34,+33)(-21,+39)(-5,+26)\put(-4,+25){\vector(1,-1){0}}
     \qbezier(-34,-33)(-21,-39)(-5,-26)\put(-4,-25){\vector(1,+1){0}}
     \qbezier(+34,+33)(+21,+39)(+5,+26)\put(35,+31){\vector(1,-1){0}}
     \qbezier(+34,-33)(+21,-39)(+5,-26)\put(35,-31){\vector(1,+1){0}}
}
\put(158,0){
\put(2,5){\vector(2,3){11}}
\put(2,-5){\vector(2,-3){11}}
\put(19,33){\qbezier(0,0)(7,8)(15,0)\put(16,-1){\vector(1,-1){0}}}
\put(19,-33){\qbezier(0,0)(7,-8)(15,0)\put(16,1){\vector(1,1){0}}}
\put(40,22){\vector(2,-3){11}}
\put(40,-22){\vector(2,3){11}}
}
}
\end{picture}
\end{minipage}\end{equation}

If exact rows or columns are introduced between the given
nonexact ones, the splicings move farther apart (as the ``staircases''
on which $C$ and $L$ lie move apart), with
a ``normal'' stretch between them.

\section{Total homology}\label{S.total_homology}
It is probably foolhardy, at very least, for someone who does not
know the theory of spectral sequences to attempt to say something
about the total homology of a double complex.
However, I shall note here some connections between
that subject and the constructions $^\bx A$ and $A_\bx$ we have
been working with.

Let us be given
\begin{equation}\begin{minipage}[c]{30pc}\label{d.complex}
a double complex with
objects $A_{i,r},$
vertical arrows $\delta_1:A_{i,r}\to A_{i+1,r},$ and
horizontal arrows $\delta_2:A_{i,r}\to A_{i,r+1},$
$(i,r\in\Z).$
\end{minipage}\end{equation}
In particular, at each object $A_{i,r},$ we have
\begin{equation}\begin{minipage}[c]{30pc}\label{d.d1d2}
$\delta_1\delta_2\ =\ \delta_2\delta_1,$\qquad
$\delta_1\delta_1\ =\ 0,$\qquad $\delta_2\delta_2\ =\ 0.$
\end{minipage}\end{equation}

At this point, one usually defines the \emph{total complex} induced
by this double complex to have for objects the direct
sums $A_n = \bigoplus_{i+r=n} A_{i,r},$ assuming the countable direct
sum construction to be defined and exact in our abelian category $\A.$
But we may as well be more general.
Let $\A^\Z$ denote the abelian category of all $\!\Z\!$-tuples
$X=(X_i)_{i\in\Z}$ of objects of $\A,$ and let $\sum:\,\A^\Z\to\B$
be \emph{any} exact functor into an abelian category $\B$ which
commutes with shift, i.e., has the property that for
$(X_i)\in\A^\Z$ we have a functorial isomorphism
\begin{equation}\begin{minipage}[c]{30pc}\label{d.shift}
$\sum_i X_i\ \cong\ \sum_i X_{i+1}.$
\end{minipage}\end{equation}

For instance, suppose $\A=\B=$ the category of all $\!R\!$-modules
for $R$ a fixed ring.
(As noted in \S\ref{S.salamander}, ``module'' can
mean either left or right module.)
Then we might take $\sum$ to be: (i)~the operator of
direct sum, or (ii)~the operator of
direct product, or (iii)~or (iv)~the right- or
left-truncated product operators, given by $\sum X_i =
(\prod_{i<0}X_i)\times(\bigoplus_{i\geq 0}X_i),$ respectively,
$\sum X_i = (\bigoplus_{i<0}X_i)\times(\prod_{i\geq 0}X_i)$
(these might be called ``formal Laurent sum'' operations; the
reader should check that they indeed satisfy~(\ref{d.shift})),
or even some very ``un-sum-like'' constructions,
such as (v)~$\sum X_i = (\prod X_i)/(\bigoplus X_i),$
or more generally (vi)~the reduced product of the $\!R\!$-modules $X_i$
$(i\in\Z)$ with respect to any translation-invariant filter on~$\Z.$

Before saying what we will do with these functors, let me digress and
point out that in any abelian category $\A$ with countable coproducts,
i.e., countable direct sums, the functor $\bigoplus_{i\in\Z}:
\A^\Z\to\A$ satisfies~(\ref{d.shift}) and is \emph{right} exact
(since it is a left adjoint, and therefore respects coequalizers);
and, dually, when the countable direct product functor $\prod_{i\in\Z}$
exists, it satisfies~(\ref{d.shift}) and is left exact.
For $\A$ the category of all $\!R\!$-modules, it is easy to check
that these constructions are exact on both sides;
but there are abelian categories in which countable products
or coproducts exist but are not exact.
For instance,
in the category $\A$ of \emph{torsion} abelian groups, infinite
products are given by the torsion subgroup of the
direct product as groups \cite[Exercise~I.8]{SL.Alg}.
It is easy to see that the direct product
in this category of the family of short exact sequences
$0\to\Z/(p^i)\to\Z/(p^{i+1})\to\Z/(p)\to 0,$ for $p$ a fixed prime
and $i$ ranging over the natural numbers,
loses exactness on the right: no element of the product of the
middle terms maps to the element $(1,1,1,\dots)$ of
the product of the right-hand terms.
(As described, this example is a family of short
exact sequences indexed by $\N;$ but
by associating copies of the zero short exact
sequence to indices $i<0,$ we get an example of non-right-exactness
of products indexed by $\Z.)$
Applying Pontryagin duality \cite[Theorem~1.7.2]{WR} to this example,
we get non-left-exactness of countable \emph{coproducts} in the
category of totally disconnected compact Hausdorff abelian groups
(though it seems harder to describe the elements involved).
Thus, these two functors, and others like them, are
\emph{excluded} as candidates for the $\sum$ we are considering.

Returning to where we left off, suppose we
are given an exact functor $\sum:\A^\Z\to\B$
satisfying~(\ref{d.shift}), and
a double complex~(\ref{d.complex}) in $\A.$
Then we define
\begin{equation}\begin{minipage}[c]{30pc}\label{d.An}
$A_n\ =\ \sum_{i\in\Z}\ A_{i,n-i},$\quad for each $n\in\Z.$
\end{minipage}\end{equation}
The families of maps $\delta_1$ and $\delta_2$ of~(\ref{d.complex})
induce maps which we shall denote by the same symbols,
\begin{equation}\begin{minipage}[c]{30pc}\label{d.d1d2total}
$\delta_1,\delta_2:\ A_n\ \to\ A_{n+1},$\quad for each $n\in\Z.$
\end{minipage}\end{equation}
(Remark: the isomorphism~(\ref{d.shift}) is needed in the definition
of $\delta_1,$ but not in that of $\delta_2;$ essentially
because we decided arbitrarily that the $i$ of the operator
$\sum_i$ in~(\ref{d.An}) would index
the \emph{first} subscript of $A_{i,n-i}.)$
These maps will again clearly satisfy~(\ref{d.d1d2}).
Since $\delta_1$ and $\delta_2$ now represent maps which
can simultaneously have the same range and the same domain
(as in~(\ref{d.d1d2total})), we can
add and subtract them, and~(\ref{d.d1d2}) immediately yields
\begin{equation}\begin{minipage}[c]{30pc}\label{d.(+)(-)}
$(\delta_1+\delta_2)\,(\delta_1-\delta_2)\ =\ 0
\ =\ (\delta_1-\delta_2)\,(\delta_1+\delta_2).$
\end{minipage}\end{equation}

Thus, if for each $n$ we let
\begin{equation}\begin{minipage}[c]{30pc}\label{d.(-1)^n}
$\delta\ =\ \delta_2+(-1)^n\delta_1:\ A_n\longrightarrow A_{n+1},$
\end{minipage}\end{equation}
we get a complex
\begin{equation}\begin{minipage}[c]{30pc}\label{d.total.cx}
\begin{picture}(210,10)
\put(0,0){$\cdots$}
\put(23,3){\vector(1,0){23}\put(-18,4){$\delta$}}
\put(50,0){$A_{n-1}$}
\put(79,3){\vector(1,0){23}\put(-18,4){$\delta$}}
\put(105,0){$A_n$}
\put(125,3){\vector(1,0){23}\put(-18,4){$\delta$}}
\put(152,0){$A_{n+1}$}
\put(181,3){\vector(1,0){23}\put(-18,4){$\delta$}}
\put(208,0){$\cdots\,.$}
\end{picture}
\end{minipage}\end{equation}
We shall call~(\ref{d.total.cx})
the \emph{total complex} (with respect to the
functor $\sum)$ of our given double complex.
Since the maps $\delta$ come from maps going downward and
to the right on our original double complex, we shall denote the
homology objects of the above complex by
\begin{equation}\begin{minipage}[c]{30pc}\label{d.diag}
$A_n\diag\ =\ \r{Ker}(A_n\stackrel{\delta}{\to} A_{n+1})\,/
\,\r{Im}(A_{n-1}\stackrel{\delta}{\to} A_n).$
\end{minipage}\end{equation}

So far, this is nothing new.
We now bring our donor and receptor objects into the picture.
Let us define
\begin{equation}\begin{minipage}[c]{30pc}\label{d.total_bx}
$A_{n\bx}\ =\sum_i\ A_{i,n-i\bx}\,,$\quad
$^\bx A_n=\ \sum_i\,{^\bx A_{i,n-i}}\,.$
\end{minipage}\end{equation}

From the exactness and shift-invariance of $\sum,$ it follows that
within each object $A_n,$ subobjects such as
$\r{Ker}(\delta_1{:A_n\to A_{n+1}})$ and
$\r{Im}(\delta_1{:A_{n-1}\to A_n})$
will be given by the corresponding ``sums'',
$\sum_i\,\r{Ker}(\delta_1{:}\ \linebreak[0] A_{i,n-i}\to A_{i+1,n-i})$
and $\sum_i\,\r{Im}(\delta_1{:A_{i-1,n-i}\to A_{i,n-i}}),$
and similarly for more complicated expressions.
One deduces that
\begin{equation}\begin{minipage}[c]{30pc}\label{d.total_bx=}
$A_{n\bx}\ =
\ \r{Ker}(\delta_1\delta_2)/(\r{Im}\,\delta_1 + \r{Im}\,\delta_2),
\quad
^\bx A_n=
\ (\r{Ker}\,\delta_1\cap\,\r{Ker}\,\delta_2)/\r{Im}(\delta_1\delta_2).$
\end{minipage}\end{equation}
(Here, in the numerators, the symbols
$\delta_1,$ $\delta_2,$ $\delta_1\delta_2,$
denote the maps so named having domain $A_n,$
and in the denominators, the maps with codomain $A_n.)$

In view of~(\ref{d.total_bx=}),
the identity morphism of $A_n$ induces intramural maps
\begin{equation}\begin{minipage}[c]{30pc}\label{d.total_intramural}
$^\bx A_n\longrightarrow A_n\diag\longrightarrow A_{n\bx}\,.$
\end{minipage}\end{equation}

(We could have defined $A_n\ver$ and $A_n\hor$ analogously
to~(\ref{d.total_bx}), noted characterizations of them
analogous to~(\ref{d.total_bx=}), and gotten a commuting diagram
\begin{equation}\begin{minipage}[c]{15pc}\label{d.<|>}
\begin{picture}(100,100)
\put(50,50){
\put(0,40){$^\bx A_n$}
\put(-40,0){$A_n\ver$}
\put(2,0){$A_n\diag$}
\put(40,0){$A_n\hor\,;$}
\put(2,-40){$A_{n\bx}$}
\put(-2,33){\vector(-1,-1){20}}
\put(13,33){\vector(0,-1){20}}
\put(28,33){\vector(1,-1){20}}
\put(-22,-7){\vector(1,-1){20}}
\put(13,-7){\vector(0,-1){20}}
\put(48,-7){\vector(-1,-1){20}}
}
\end{picture}
\end{minipage}\end{equation}
but we shall not need these additional objects.)

Finally, the two sets of \emph{extra}mural maps constructed
from our original double complex in \S\ref{S.salamander},
combined with~(\ref{d.total_bx}), yield,
for each $n,$ two maps which we shall call
\begin{equation}\begin{minipage}[c]{30pc}\label{d.bd12}
\begin{picture}(140,20)
\put(0,4){$A_{n\bx}$}
\put(27,9){\vector(1,0){30}\put(-20,4){$\bd_1$}}
\put(27,4){\vector(1,0){30}\put(-20,-13){$\bd_2$}}
\put(62,4){$^\bx A_{n+1}.$}
\end{picture}
\vspace{.8em}
\end{minipage}\end{equation}
In terms of the description~(\ref{d.total_bx=}) of
$A_{n\bx}$ and $^\bx A_{n+1},$ we see that $\bd_1$ and $\bd_2$
are induced in $\B$ by $\delta_1,\ \delta_2:A_n\to A_{n+1}.$

Let us now write (analogous to~(\ref{d.(-1)^n})),
\begin{equation}\begin{minipage}[c]{30pc}\label{d.bd}
$\bd\ =\ \bd_2+(-1)^n\,\bd_1:\ A_{n\bx}\longrightarrow{^\bx A_{n+1}}.$
\end{minipage}\end{equation}
We find that the composite
of this map with the first intramural map of~(\ref{d.total_intramural}),
\begin{equation}\begin{minipage}[c]{30pc}\label{d.bd.intermur}
\begin{picture}(140,20)
\put(0,4){$A_{n-1\bx}$}
\put(37,6){\vector(1,0){30}\put(-20,4){$\bd$}}
\put(70,4){$^\bx A_n$}
\put(93,6){\vector(1,0){30}}
\put(125,2){$A_n\diag$}
\end{picture}
\end{minipage}\end{equation}
is zero, since the first arrow maps into the denominator
of~(\ref{d.diag}).
This says that the two composites
\begin{equation}\begin{minipage}[c]{30pc}\label{d.bd12.intermur}
\begin{picture}(140,20)
\put(0,4){$A_{n-1\bx}$}
\put(37,9){\vector(1,0){30}\put(-20,4){$\bd_1$}}
\put(37,4){\vector(1,0){30}\put(-20,-13){$\bd_2$}}
\put(70,4){$^\bx A_n$}
\put(93,6){\vector(1,0){30}}
\put(125,4){$A_n\diag$}
\end{picture}
\vspace{.8em}
\end{minipage}\end{equation}
agree up to sign.
Hence, below, we shall just refer to the composite involving $\bd_1.$
The same comments apply to the composites
\begin{equation}\begin{minipage}[c]{30pc}\label{d.intermur.bd12}
\begin{picture}(140,20)
\put(0,4){$A_n\diag$}
\put(25,6){\vector(1,0){30}}
\put(58,4){$A_{n\bx}$}
\put(83,9){\vector(1,0){30}\put(-20,4){$\bd_1$}}
\put(83,4){\vector(1,0){30}\put(-20,-13){$\bd_2$}}
\put(115,4){$^\bx A_{n+1}.$}
\end{picture}
\vspace{1em}
\end{minipage}\end{equation}

We can now state a version of the Salamander Lemma
for total homology.

\begin{lemma}\label{L.total_exact}
In the above context, for each $n$ the
$\!6\!$-term sequence of intramural and
extramural maps and their composites
\begin{equation}\begin{minipage}[c]{35pc}\label{d.total_exact}
\begin{picture}(370,27)
\put(010,10){$A_{n-1\bx}$}
\put(045,13){\vector(1,0){55}}
\put(050,18){$\bd_1$}
\put(065,18){$\scriptstyle{^\bx A_n}$}
\put(060,3){\em\scriptsize cf.~(\ref{d.bd12.intermur})}
\put(070,13){\middot}
\put(103,10){$A_n\diag$}
\put(127,13){\vector(1,0){30}}
\put(160,10){$A_{n\bx}$}
\put(185,13){\vector(1,0){25}}
\put(213,10){$^\bx A_{n+1}$}
\put(247,13){\vector(1,0){25}}
\put(275,10){$A_{n+1}\diag$}
\put(308,13){\vector(1,0){55}}
\put(313,18){$\scriptstyle{A_{n+1\bx}}$}
\put(345,18){$\bd_1$}
\put(323,3){\em\scriptsize cf.~(\ref{d.intermur.bd12})}
\put(333,13){\middot}
\put(366,10){$^\bx A_{n+2}$}
\end{picture}
\end{minipage}\end{equation}
is exact.
\end{lemma}

\proof
Rather than re-proving this, let us deduce it from the Salamander
Lemma (Lemma~\ref{L.salamander}) by a trick.
We define a double complex $(B_{i,r})$ in $\B$ in which
$B_{i,r}=A_{i+r},$ the horizontal maps are the maps $\delta,$
and the vertical maps are the $\delta_1.$
Thus, letting $n=i+r,$ we have
\begin{equation}\begin{minipage}[c]{35pc}\label{d.Bir=Ai+r}
\begin{picture}(460,110)
\put(40,100){$B_{i-1,r}$}
\put(40,50){$B_{i,r}$}
\put(90,50){$B_{i,r+1}$}
\put(90,0){$B_{i+1,r+1}$}
\multiput(12,53)(50,0){2}{\vector(1,0){26}\put(-20,3){$\delta$}}
\put(122,53){\vector(1,0){26}\put(-20,3){$\delta$}}

\multiput(48,43)(50,0){2}{\multiput(0,0)(0,50){2}{%
	\vector(0,-1){30}\put(1,-15){$\delta_1$}}}

\multiput(258,43)(60,0){2}{\multiput(0,0)(0,50){2}{%
	\vector(0,-1){30}\put(1,-15){$\delta_1$}}}

\put(250,100){$A_{n-1}$}
\put(250,50){$A_n$}
\put(310,50){$A_{n+1}$}
\put(310,0){$A_{n+2}$}
\put(210,53){\vector(1,0){38}\put(-37,5){$\delta_2\mp\delta_1$}}
\put(270,53){\vector(1,0){38}\put(-37,5){$\delta_2\pm\delta_1$}}
\put(340,53){\vector(1,0){38}\put(-37,5){$\delta_2\mp\delta_1$}}

\multiput(175,48)(0,3){2}{\line(1,0){10}}
\end{picture}
\end{minipage}\end{equation}

From~(\ref{d.d1d2}) and~(\ref{d.(-1)^n}), we see that
$\delta_1\delta=\delta_1\delta_2=\delta_2\delta_1=\delta\delta_1,$
and that in each object,
$\r{Im}\,\delta_1+\r{Im}\,\delta=\r{Im}\,\delta_1+\r{Im}\,\delta_2$ and
$\r{Ker}\,\delta_1\cap\r{Ker}\,\delta=\r{Ker}\,\delta_1\cap
\r{Ker}\,\delta_2.$
It easily follows that $B_{i,r}\hor=A_n\diag,$
$B_{i,r\bx}=A_{n\bx},$ $^\bx B_{i,r}={^\bx A_n},$
and that~(\ref{d.total_exact}) is just the $\!6\!$-term
exact sequence associated with the above horizontal arrow of this
double complex.
\endproof

Likewise, Corollary~\ref{C.1-row-long}, applied to any row
of the complex~(\ref{d.Bir=Ai+r}) gives:

\begin{lemma}\label{L.total_long}
If our original double complex has exact columns, then
one has a long exact sequence in $\B$:
\begin{equation}\begin{minipage}[c]{30pc}\label{d.total_long}
\begin{picture}(210,15)
\put(0,0){$\cdots$}
\put(20,3){\vector(1,0){25}\put(-18,4){$\bd$}}
\put(50,0){$^\bx A_n$}
\put(75,3){\vector(1,0){25}}
\put(105,0){$A_n\diag$}
\put(130,3){\vector(1,0){25}}
\put(160,0){$A_{n\bx}$}
\put(186,3){\vector(1,0){25}\put(-18,4){$\bd$}}
\put(212,0){$^\bx A_{n+1}$}
\put(248,3){\vector(1,0){25}}
\put(280,0){$\cdots\,.$}
\vspace{1em}
\end{picture}
\endproof
\end{minipage}\end{equation}
\end{lemma}

(Note the curious property of this sequence: that the objects
joined by the connecting morphisms $\bd$ are \emph{isomorphic}
under a different map, $\bd_1,$ by Corollary~\ref{C.extra-iso}.)

\begin{corollary}\label{C.diag0}
If our original double complex has \emph{exact rows and columns}, and
is \emph{weakly bounded} \textup{(}e.g., if $A_{i,r}=0$ whenever $i$
or $r$ is negative\textup{)}, then all total homology objects
$A_n\diag$ are zero.
\end{corollary}

\proof
In the original double complex $(A_{i,r}),$ the donor
and receptor objects form isomorphic chains
going upward to the right and downward to the left
by Corollary~\ref{C.extra-iso} and our exactness assumptions;
hence, by weak boundedness, they are all zero.
Thus, for all~$n,$ $^\bx A_n = \sum {^\bx A_{i,n-i}}=\sum 0=0,$
and similarly $A_{n\bx}=0.$
So by~(\ref{d.total_long}), $A_n\diag=0.$
\endproof

There is no more to be said about total homology under this hypothesis,
so to finish this section, let us return to the weaker hypothesis
of Lemma~\ref{L.total_long}, and examine the behavior
of the exact sequence~(\ref{d.total_long}) for various
choices of $\A$ and $\sum.$

\begin{convention}\label{Cn.exact_cols}
For the remainder of this section, we shall assume as
in Lemma~\ref{L.total_long} that our given double
complex $(A_{i,r})_{i,r\in\Z}$ has exact columns.
\end{convention}

Consider now the pair of objects $A_{n\bx}$ and $^\bx A_{n+1}$
of the exact sequence~(\ref{d.total_long}).
We have two maps, $\bd_1$ and $\bd_2,$
between them, as in~(\ref{d.bd12}),
of which $\bd_1$ is now an isomorphism because of our assumption
of exact columns (Corollary~\ref{C.extra-iso}).
Let us think of the objects $A_{n\bx}$ and $^\bx A_{n+1}$
as the combinations, under the functor $\sum,$ of the
donor objects, respectively the receptor objects, in
the chain of objects and maps in $\A$
lying on one of the ``staircases'' in~(\ref{d.linked_long}).
Flattened out, such a staircase looks like
\begin{equation}\begin{minipage}[c]{32pc}\label{d.bd1bd2}
\begin{picture}(340,30)
\put(0,10){$\cdots$}
\put(20,13){\vector(1,0){25}\put(-18,4){$\bd_2$}}
\put(50,10){$^\bx A_{i+1,n-i}\,\cong$\put(-10,13){$\bd_1^{-1}$}}
\put(115,10){$A_{i,n-i\bx}$}
\put(155,13){\vector(1,0){25}\put(-18,4){$\bd_2$}}
\put(185,10){$^\bx A_{i,n-i+1}\,\cong$\put(-10,13){$\bd_1^{-1}$}}
\put(250,10){$A_{i-1,n-i+1\bx}$}
\put(310,13){\vector(1,0){25}\put(-18,4){$\bd_2$}}
\put(340,10){$\cdots\,.$}
\vspace{1em}
\end{picture}
\end{minipage}\end{equation}
Here the maps $\bd_1,\bd_2: A_{n\bx}\to {^\bx A_{n+1}}$
between our ``total'' donor and receptor objects
arise under $\sum$ from the above maps of~(\ref{d.bd1bd2}).

If we compose each arrow $\bd_2$ of~(\ref{d.bd1bd2}) with the
preceding (or the
following) isomorphism $\bd_1^{-1},$ we may regard~(\ref{d.bd1bd2})
as a directed system in $\A.$
This suggests that we look at the direct limit
$\injLim_{i\to-\infty}\,A_{i,n-i\bx}\ =
\ \injLim_{i\to-\infty}\,{^\bx A_{i,n-i+1}}$ of that system,
if this exists in $\A.$

If in fact $\A$ has countable coproducts (direct sums),
then it has general countable colimits
(by \cite[Proposition~7.6.6]{245}, and the fact
that, being an abelian category, it has coequalizers)
and so, in particular, countable direct limits.
Examining how this is proved, we see that in this
situation, the direct limit
of~(\ref{d.bd1bd2}) can be constructed as the cokernel of the map
\begin{equation}\begin{minipage}[c]{30pc}\label{d.oplus_map}
$\bd_1-\bd_2:
\ \bigoplus_i\,A_{i,n-i\bx}\ \longrightarrow
\ \bigoplus_i\,{^\bx A_{i,n-i+1}}.$
\end{minipage}\end{equation}

Now suppose countable coproducts are exact in $\A,$
and take $\B=\A$ and $\sum=\bigoplus.$
Then (up to a possible change of sign in $\bd_1,$ which clearly
won't change the direct limit of~(\ref{d.bd1bd2})) we see
that~(\ref{d.oplus_map}) is just
\begin{equation}\begin{minipage}[c]{30pc}\label{d.just_bd}
$\bd:\ A_{n\bx}\ \longrightarrow\ ^\bx A_{n+1}.$
\end{minipage}\end{equation}

In summary: if $\sum=\bigoplus,$ then the cokernel
of the step~(\ref{d.just_bd}) in the exact
sequence~(\ref{d.total_long})
is given by the direct limit of~(\ref{d.bd1bd2}).

The \emph{kernel} of~(\ref{d.just_bd}) does not in
general have such a natural description for $\sum=\bigoplus.$
But if $\A$ is the category of $\!R\!$-modules
$(R$ any ring), that kernel will be zero!
Indeed, consider any nonzero $x\in A_{n\bx}=\bigoplus_i\,A_{i,n-i\bx}.$
Let $i$ be the largest integer such that $0\neq x_i\in A_{i,n-i\bx}$
(corresponding diagrammatically to the \emph{lowest} position
where $A_{n\bx}$ has nonzero component on the staircase where it lives).
Examining the $^\bx A_{i+1,n-i}$ component of
$\bd(x),$ we see that this is $\bd_1(x_i),$ which is nonzero
because $\bd_1$ is an isomorphism.
So $\bd(x)$ is injective, as claimed.

Applying these observations to the exact sequence~(\ref{d.total_long}),
we get:

\begin{corollary}\label{C.total_mdls_sum}
If $R$ is a ring, and $(A_{i,r})$ a double complex of $\!R\!$-modules
with exact columns, then its total homology
with respect to $\bigoplus$ is described by
\begin{equation}\begin{minipage}[c]{33pc}\label{d.total_mdls}
$A_n\diag\ \cong
\ \r{Cok}(\bd:A_{n-1\bx}\to{^\bx A_n})\ \cong
\ \injLim_{i\to-\infty}\,A_{i,n-i-1\bx}\ \cong
\ \injLim_{i\to-\infty}\,{^\bx A_{i,n-i}},$
\end{minipage}\end{equation}
where the direct limits
are over the system~\textup{(\ref{d.bd1bd2})}.\endproof
\end{corollary}

Dualizing the observation following~(\ref{d.just_bd}),
we see that if $\A$ an abelian category with countable direct products
and these are exact, and we take $\sum=\prod,$
then the kernel of~(\ref{d.just_bd}) is the \emph{inverse}
limit of~(\ref{d.bd1bd2}).
The \emph{cokernel} is now hard to describe, even
for $\A$ the category of $\!R\!$-modules; but
I shall develop the description in that case
(with $n-1$ in place of $n,$ for convenience) below.

We begin by noting that for general $\A$ and $\sum,$ the exactness
of~(\ref{d.total_long}) at $^\bx A_n$ and
$A_n\diag$ tells us that
\begin{equation}\begin{minipage}[c]{30pc}\label{d.3ways}
$\r{Cok}(\bd:A_{{n-1}\bx}\to{^\bx A_n})\ \cong
\ \r{Im}(^\bx A_n\to A_n\diag)\ =
\ \r{Ker}(A_n\diag\to A_{n\bx}).$
\end{minipage}\end{equation}
(where the arrows in the last two expressions are intramural maps).

Now let $\A$ again be the category of $\!R\!$-modules,
let $(A_{i,r})$ still be a double complex in $\A$ with exact columns,
and let its total complex, its total homology,
and our related objects
be defined using the functor $\sum=\prod$ on $\A.$
Suppose an element $x\in A_n\diag$ has the property
that for every \emph{finite} subset $I\subseteq\Z,$ $x$ can be
represented by a cycle in $A_n$ which has zero component in
$A_{i,n-i}$ for all $i\in I.$
Then I will call $x$ a ``peekaboo element'' (because wherever
you look, it isn't there!)
The set of these elements forms a submodule of $A_n\diag,$
which I shall denote $\r{PB}(A_n\diag).$
I claim that~(\ref{d.3ways}), as a submodule of $A_n\diag,$
is precisely $\r{PB}(A_n\diag).$

Let us first show
that $\r{PB}(A_n\diag)$ is contained in~(\ref{d.3ways}),
looked at as $\r{Ker}(A_n\diag\to A_{n\bx}).$
Suppose $x\in\r{PB}(A_n\diag),$ and let $x$ be
represented by a cycle
\begin{equation}\begin{minipage}[c]{30pc}\label{d.x_i}
$(x_i)\in\prod_i A_{i,n-i}.$
\end{minipage}\end{equation}
Since $x$ is a peekaboo element, for any $j\in\Z$ we can
modify $(x_i)$ by a boundary $\delta(y_i)$ to get an
element having $\!j\!$-th coordinate $0.$
But doing this changes $x_j$ by $\delta_2(y_j)\pm\delta_1(y_{j-1}),$
so if it can send $x_j$ to zero,
we must have $x_j\in \delta_2(A_{j,n-j-1})+\delta_1(A_{j-1,n-j}),$
which means that $x_j$ has zero image in $A_{j,n-j\bx}\,.$
Since we have proved this for arbitrary $j,$ the
image of $x$ in $A_{n\bx}=\prod A_{i,n-i\bx}$ is zero,
i.e., $x\in\r{Ker}(A_n\diag\to A_{n\bx}),$ as claimed.

Conversely, suppose $x\in A_n\diag$ lies in~(\ref{d.3ways}),
which we now look at as $\r{Im}(^\bx A_n\to A_n\diag).$
This says that $x$ can be represented by a cycle $(x_i)$ as
in~(\ref{d.x_i}) such that each coordinate $x_i$
is annihilated by both $\delta_1$ and $\delta_2.$
We wish to show that for any finite subset $I\subseteq\Z,$
our cycle $(x_i)$ can be modified by a
boundary so that it becomes zero at each coordinate in $I.$
Suppose inductively that we have been able to achieve zero entries at
all coordinates in $I-\{j\},$ where $j=\min(I)$ (corresponding to the
\emph{highest} point we are interested in on our
staircase), in the process preserving the condition that each coordinate
is annihilated by $\delta_1$ and $\delta_2.$
For notational simplicity, let us again call this modified element
$(x_i).$
In particular, the coordinate $x_j$ is annihilated by $\delta_1;$
so by exactness of the columns of our double complex, we can
write $x_j=\delta_1(z_{j-1})$ for some $z_{j-1}\in A_{j-1,n-j}.$
If we let $(z_i)\in A_{n-1}$ be the element with $\!j{-}1\!$-st
coordinate $z_{j-1}$ and all other coordinates $0,$ it is
easy to verify that $(x_i)+(-1)^n\delta((z_i))$ has zero coordinates
at all indices in $I.$
(The two coordinates that differ from those of $(x_i)$ are
the $\!j\!$-th, which has been brought to zero, and
the $\!j{-}1\!$-st, which we don't care about because $j-1\notin I.)$
To see that all coordinates of this element are still
annihilated by $\delta_1$ and $\delta_2,$ note that it
suffices to prove that $\delta((z_i))$ has this property, for which,
by~(\ref{d.d1d2}) and~(\ref{d.(-1)^n}),
it will suffice to show that $\delta_2\delta_1(z_{j-1})=0.$
But $\delta_2\delta_1(z_{j-1})=\delta_2(x_j)=0$ by choice of $(x_i).$
This completes the proof that~(\ref{d.3ways})
is given by $\r{PB}(A_n\diag).$

Inserting into~(\ref{d.total_long}) our earlier description of
the kernel of $\bd,$ and this description of its cokernel, we get

\begin{corollary}\label{C.total_mdls_prod}
If $R$ is a ring, and $(A_{i,r})$ a double complex of $\!R\!$-modules
with exact columns, and we form its total homology objects $A_n\diag$
with respect to $\prod,$ then for each $n$ we have a short exact
sequence
\begin{equation}\begin{minipage}[c]{35pc}\label{d.PB}
$0\ \longrightarrow \ \r{PB}(A_n\diag)
\ \longrightarrow\ A_n\diag\ \longrightarrow
\ (\projLim_{i\to\infty}\,A_{i,n-i\bx}
\,\cong\,\projLim_{i\to\infty}{^\bx A_{i,n-i+1}})\ \longrightarrow\ 0,$
\end{minipage}\end{equation}
where the inverse limits
are over the system~\textup{(\ref{d.bd1bd2})}.\endproof
\end{corollary}

For an example in which the term $\r{PB}(A_n\diag)$ of~(\ref{d.PB})
is nonzero, let $p$ be any prime, and
consider the double complex of abelian groups,
\begin{equation}\begin{minipage}[c]{15pc}\label{d.PB_2cx}
\begin{picture}(150,150)
\multiput(0,0)(30,30){4}{
	\put(30,54){\vector(0,-1){17}}
	\put(6,30){\vector(1,0){18}\put(1,-3){$\Z$}}
	\put(36,30){\vector(1,0){18}\put(-14,4){$p$}}
	\put(30,8.5){\multiput(0,0)(2,0){2}{\line(0,1){15}}}
	}
\multiput(0,0)(30,30){3}{
	\put(30,84){\vector(0,-1){17}}
	\put(6,60){\vector(1,0){18}\put(2,-3){$0$}}
	\put(55,27){$\Z$}
	\put(60,24){\vector(0,-1){17}}
	\put(66,30){\vector(1,0){18}}
	}
\multiput(0,0)(30,30){2}{
	\put(30,114){\vector(0,-1){17}}
	\put(6,90){\vector(1,0){18}\put(2,-3){$0$}}
	\put(86,27){$0$}
	\put(90,24){\vector(0,-1){17}}
	\put(96,30){\vector(1,0){18}}
	}
	\put(30,144){\vector(0,-1){17}}
	\put(6,120){\vector(1,0){18}\put(2,-3){$0$}}
	\put(116,27){$0$}
	\put(120,24){\vector(0,-1){17}}
	\put(126,30){\vector(1,0){18}}
\put(150,10){,}
\end{picture}
\end{minipage}\end{equation}
where the arrows labeled ``$p$'' represent multiplication
by $p,$ and the vertical equals-signs denote the identity map.

Let $n$ be the common value of the sum of the subscripts on
the objects $\Z$ in
the lower diagonal string of such objects in~(\ref{d.PB_2cx}),
so that $A_n$ is the direct product of these objects.
Then $A_{n+1}=0,$ so all members of $A_n$ are cycles.
For the same reason, the second inverse limit in~(\ref{d.PB})
is zero; so~(\ref{d.PB}) shows that
all elements of $A_n\diag$ are peekaboo elements.
To establish the existence of \emph{nonzero}
peekaboo elements, we must therefore show that $A_n\diag\neq 0.$

To this end, consider the map $\sigma$ from $A_n$ to the
group $\Z_p$ of \emph{$\!p\!$-adic integers}, given by
\begin{equation}\begin{minipage}[c]{30pc}\label{d.to_p-adic}
$\sigma((x_i)_{i\in\Z})\ =\ \sum^{\infty}_{i=0}\ (-1)^{in}\,p^i\,x_i
\ \in\ \Z_p.$
\end{minipage}\end{equation}
(Note the ``cut-off'': all $x_i$ with $i<0$ are ignored.)

If $(x_i)\in A_n$ is a boundary, $(x_i)=\delta((y_i)),$ then
for each $i$ we have
$x_i=\delta_2(y_i)+(-1)^{n-1}\delta_1(y_{i-1})=p\,y_i-(-1)^n\,y_{i-1}.$
This makes the computation~(\ref{d.to_p-adic}) of $\sigma((x_i))$
a ``telescoping sum'', where all terms cancel except $-(-1)^n y_{-1}.$
Thus, for every boundary $(x_i)=\delta((y_i)),$ the element
$\sigma((x_i))\in\Z_p$ belongs to $\Z.$
On the other hand, we can clearly choose elements
(hence, cycles) $(x_i)\in A_n$
for which $\sigma((x_i))$ is an arbitrary member of $\Z_p.$
Thus, there are cycles in $A_n$ which are not boundaries, hence
nonzero peekaboo elements in $A_n\diag.$

It is not hard to show, conversely, that every $(x_i)\in A_n$
such that $\sigma(x_i)\in\Z$ is a boundary, and to conclude that
$A_n\diag\cong\Z_p/\Z.$

(Incidentally, in~(\ref{d.PB_2cx}) we could replace $\Z$ by
a polynomial ring $k[x],$ for $k$ a field, and the element $p$ by $x.$
In place of the $\Z_p/\Z$ in the final result we would then
get $k[[x]]/k[x].$
Regarding this example as a double complex of
\emph{$\!k\!$-vector-spaces,} we see that the existence of nonzero
peekaboo elements in $A_n\diag$ does not require any
noncompleteness property of the base ring.)

We now come to our last choice of $\sum$ on the category $\A$ of
$\!R\!$-modules, the right truncated product (or ``left formal
Laurent sum'') functor.
Thus we let
\begin{equation}\begin{minipage}[c]{30pc}\label{d.trunc_A_n}
$A_n\ =\ (\prod_{i<0}\,A_{i,n-i})\ \times
\ (\bigoplus_{i\geq 0}\,A_{i,n-i}).$
\end{minipage}\end{equation}
(This was example~(iii) in the sentence following~(\ref{d.shift}).
I use ``right'' with reference to the subscript $i$ in
$\sum X_i=(\prod_{i<0}\,X_i)\times(\bigoplus_{i\geq 0}\,X_i),$ thinking
of the indices as written in increasing order from left to right.
Unfortunately, since in our double complex the first subscript indexes
the row, and increases in the downward direction in our
diagrams, ``\emph{right} truncated'' means, for these diagrams,
``truncated in the \emph{downward left} direction along
each diagonal''.)
We shall find that this yields the simplest, i.e., most trivial,
characterization of the total homology.

Note that if in $\bd=\bd_2\pm\bd_1: A_{n-1\bx}\to{^\bx A_n},$
we look at the effects of the respective operators $\bd_1$ and
$\bd_2$ on the first subscript of $A_{i,n-1-i\bx}\,,$
each operator adds to this subscript a constant, namely
$1,$ respectively $0.$
We can think of this as saying our operators are
each ``homogeneous''; but they are homogeneous of distinct degrees,
and the operator of higher degree, $\bd_1,$ is invertible.
It follows that $\bd$ will be invertible!
Indeed, let us write $\bd=(1-\varepsilon)(\pm\bd_1),$ where
$\varepsilon=\pm\bd_2\bd_1^{-1},$ noting that $\varepsilon$
is homogeneous of degree $-1.$
Then we see that the formal inverse
\begin{equation}\begin{minipage}[c]{30pc}\label{d.bd-1}
$\bd^{-1}\,=\ \pm\,\bd_1^{-1}(1+\varepsilon+\varepsilon^2+\dots)$
\end{minipage}\end{equation}
converges on our right truncated product modules~(\ref{d.trunc_A_n}),
and thus gives a genuine inverse to $\bd.$
Lemma~\ref{L.total_long} now immediately gives:

\begin{corollary}\label{C.trunc_triv}
Every double complex of $\!R\!$-modules with exact
columns has trivial total homology with respect to the right
truncated product \textup{(}left formal Laurent sum\textup{)}
functor $\prod_{i<0}\times\bigoplus_{i\geq 0}.$\endproof
\end{corollary}

For the \emph{left} truncated product
functor $\bigoplus_{i<0}\times\prod_{i\geq 0}$
(the right formal Laurent sum), the kernel and cokernel of
$\delta$ seem much more difficult to describe in terms of the directed
system~(\ref{d.bd1bd2}), and I will not try to do so.
(Of course, if, reversing Convention~\ref{Cn.exact_cols},
we take rows rather than columns
exact in our double complex of $\!R\!$-modules, the behaviors
of left and right truncated products are reversed.)

In a different direction, one finds that for each of the above four
constructions $\sum,$ on \emph{any} abelian category $\A$
where it is defined and exact,
\emph{weakly bounded} double complexes with exact
columns always have trivial total homology:

\begin{corollary}\label{C.wk_bdd}
Let $\A$ be any abelian category with countable
coproducts, respectively countable products, respectively both,
which are exact functors; let $(A_{i,r})$ be a
weakly bounded double complex in $\A$ with exact columns, and suppose
we form its total complex with respect to the
coproduct functor, respectively, the product functor,
respectively, the right- or left-truncated product functor.
Then in each case, the induced maps $\bd:A_{n\bx}\to{^\bx A_{n+1}}$
\textup{(}see~\textup{(\ref{d.bd}))} will be isomorphisms,
and hence the total homology will be zero.
\end{corollary}

\subsubsection*{\textsc{Sketch of Proof}}
The weak boundedness hypothesis has the effect that for
each $n,$ the string of objects and maps
\begin{equation}\begin{minipage}[c]{32pc}\label{d.bd1bd2<->}
\begin{picture}(400,30)
\put(0,10){$\cdots$}
\put(20,13){\vector(1,0){25}\put(-18,4){$\bd_2$}}
\put(50,10){$^\bx A_{i+1,n-i}$}
\put(130,13){\vector(-1,0){30}
	\put(10,4){$\bd_1$}\put(10,-10){$\cong$}}
\put(135,10){$A_{i,n-i\bx}$}
\put(175,13){\vector(1,0){25}\put(-18,4){$\bd_2$}}
\put(205,10){$^\bx A_{i,n-i+1}$}
\put(285,13){\vector(-1,0){30}
	\put(10,4){$\bd_1$}\put(10,-10){$\cong$}}
\put(290,10){$A_{i-1,n-i+1\bx}$}
\put(350,13){\vector(1,0){25}\put(-18,4){$\bd_2$}}
\put(380,10){$\cdots$}
\vspace{1em}
\end{picture}
\end{minipage}\end{equation}
breaks up into (generally infinitely many) finite substrings,
separated by zero objects.
Thus, the map $\bd:A_{n\bx}\to{^\bx A_{n+1}}$ becomes a
$\!\sum\!$-sum of maps from finite direct sums of consecutive
donor objects to finite direct sums of consecutive receptor objects.
(One verifies this by examining how each of our four functors
$\sum$ behaves
with respect to decompositions of its domain into finite subfamilies.)
On each of these pairs of finite sums, one verifies that
the restriction of $\bd$ is invertible, by
a finite version of the computation~(\ref{d.bd-1}).
A $\!\sum\!$-sum of invertible maps is invertible,
completing the proof.
\endproof

I have not investigated the behavior of total homology
with respect to any other functors $\sum.$
In particular, I do not know of any examples of
exact functors $\sum:\A^\Z\to\B$ satisfying~(\ref{d.shift})
for which the analog of Corollary~\ref{C.wk_bdd} fails.
Cf.\ Corollary~\ref{C.diag0}, which says that any such
functor does give trivial homology on double complexes
having \emph{both} exact rows and exact columns.

\section{Further notes}\label{S.remarks}

\subsection{A formally simpler approach}\label{SS.1-obj}
A more sophisticated formulation of the basic ideas of this note
(say, of Definition~\ref{D.4things} through Lemma~\ref{L.salamander},
plus Corollary~\ref{C.linked_long} and
Lemmas~\ref{L.total_exact} and~\ref{L.total_long})
would refer to a single object $A$ of an abelian category,
possibly graded, with two commuting (or
anticommuting) square-zero endomorphisms $\delta_1$ and $\delta_2.$
The situation we have been studying would be the particular case where
the category is the bigraded additive category of double complexes
in our $\A.$
(Bigraded because we would allow subscript-shifting as well
as subscript-preserving morphisms.)
For instance, in the formulation of Corollary~\ref{C.linked_long},
the vertical exactness assumption would simply take the
form $\r{Im}(\delta_1)=\r{Ker}(\delta_1),$ and the diagram
obtained,~(\ref{d.total_long}), would reduce to what is
called in \cite{H+S} an exact couple:
\begin{equation}\begin{minipage}[c]{10pc}\label{d.exact_triangle}
\begin{picture}(50,50)
\put(-2,35){$A_\bx$}
\put(15,39){\vector(1,0){30}}
\put(47,35){$^\bx A$}
\put(48,32){\vector(-2,-3){14}}
\put(22,11){\vector(-2,3){12}}
\put(22,0){$A\,{\diag}\ ;$}
\end{picture}
\vspace{1em}
\end{minipage}\end{equation}
alongside which we would have an isomorphism $^\bx A\cong A_\bx.$
But I will not attempt to develop any of the material in this form.

\subsection{Triple complexes}\label{SS.3-cx}
At every object $A$ of an ordinary complex, one has a single
homology object, while above we have associated to each $A$ in
a \emph{double} complex four
homology objects, $A_\bx\,,$ $A\hor,$ $A\ver$ and $^\bx A.$
What would be the analogous constructions in a triple complex?

To answer this, let us examine how the four constructions
we associate to a double complex arise.
Let us start with a picture \raisebox{-.80em}{\twosquare{}\ },
representing an object $A$ of our double complex,
together with the three other objects of the
complex from which the double complex structure gives us
possibly nontrivial maps into $A,$ and the three into which
it gives us possibly nontrivial maps from $A$ (cf.~(\ref{d.fig8}));
and let us mark with dots those of these objects
occurring in the definition of each of our homology objects:
\begin{equation}\begin{minipage}[c]{30pc}\label{d.2cx_outline}
\begin{picture}(400,120)
\put(200,10){
\put(20,95){Diagrams}
\put(30,60){\twosquare{
\put(0,20){\bigdot}
\put(10,0){\bigdot}
\put(20,10){\bigdot}}}

\put(0,30){\twosquare{
\put(0,10){\bigdot}
\put(20,10){\bigdot}}}

\put(60,30){\twosquare{
\put(10,20){\bigdot}
\put(10,0){\bigdot}}}

\put(30,0){\twosquare{
\put(20,0){\bigdot}
\put(10,20){\bigdot}
\put(0,10){\bigdot}}}
}

\put(35,15){
\put(0,90){Homology objects}
\put(33,70){$^\bx A$}
\put(35,65){\vector(-1,-1){18}}
\put(52,65){\vector(1,-1){18}}
\put(3,35){$A\hor$}
\put(69,35){$A\ver$}
\put(70,30){\vector(-1,-1){18}}
\put(17,30){\vector(1,-1){18}}
\put(36,0){$A_\bx$}
}
\end{picture}
\end{minipage}\end{equation}

For each of the above diagrams, the associated homology
object is the quotient of the intersection of the kernels
of the maps from $A$ to the marked objects in the lower square,
by the sum of the images in $A$ of the marked objects in
the upper square.
The dots in the two squares are in each case located so that the maps
from the marked objects of the upper square into the marked objects
in the lower square via $A$ are all zero, i.e., so that the
image in question is indeed contained in the kernel named.
On the other hand, the dots in the upper square are in each case as low
and as far to the right as they can be without violating this
condition, given the positions of the dots in the lower square, and
the dots in the lower square are as high and as far to the left
as they can be, given the positions of the dots in the upper square.

In fact, if we partially order the vertices of our diagrams
by considering the arrows of our double
complex to go from higher to lower elements, and we supplement
our dots in the lower square with all the dots below them under
this ordering, and
those in the upper square with all the dots above them (noting that
so doing does not change the resulting sums of images, or
intersections of kernels),
\begin{equation}\begin{minipage}[c]{30pc}\label{d.2cx_redundant}
\begin{picture}(400,90)
\put(200,10){
\put(30,60){\twosquare{
\put(0,20){\bigdot}
\put(10,0){\bigdot}
\put(20,0){\bigdot}
\put(20,10){\bigdot}}}

\put(0,30){\twosquare{
\put(0,10){\bigdot}
\put(0,20){\bigdot}
\put(20,0){\bigdot}
\put(20,10){\bigdot}}}

\put(60,30){\twosquare{
\put(0,20){\bigdot}
\put(10,20){\bigdot}
\put(20,0){\bigdot}
\put(10,0){\bigdot}}}

\put(30,0){\twosquare{
\put(20,0){\bigdot}
\put(10,20){\bigdot}
\put(0,20){\bigdot}
\put(0,10){\bigdot}}}
}
\end{picture}
\end{minipage}\end{equation}
then we see that the arrays of dots in the lower squares of the
above four diagrams are precisely the four proper nonempty down-sets
(sets closed under $\leq)$ of that partially ordered set, and the
arrays of dots in the upper squares are (if we momentarily superimpose
the upper and lower squares) the \emph{complementary} up-sets.

Knowing this, we can see what the analog should be for triple complexes.
One finds that the set of vertices of a cube has $18$ proper nonempty
down-sets, each with its complementary up-set.
Under permutation of the three coordinates, these $18$
complementary pairs form $8$ equivalence classes.
For simplicity, we show on the right below only one representative
of each equivalence class, and we again show by dots only maximal
elements of each down-set and minimal elements of each up-set,

since these correspond to the objects actually needed to compute
our homology objects.
We mark each diagram $\times 1$ or $\times 3,$ to indicate the
size of its orbit under permutations of coordinates.
On the
left, we show the full partially ordered system of these objects.
The lines showing the order-relations in that partially
ordered set induce intramural maps among our homology objects.
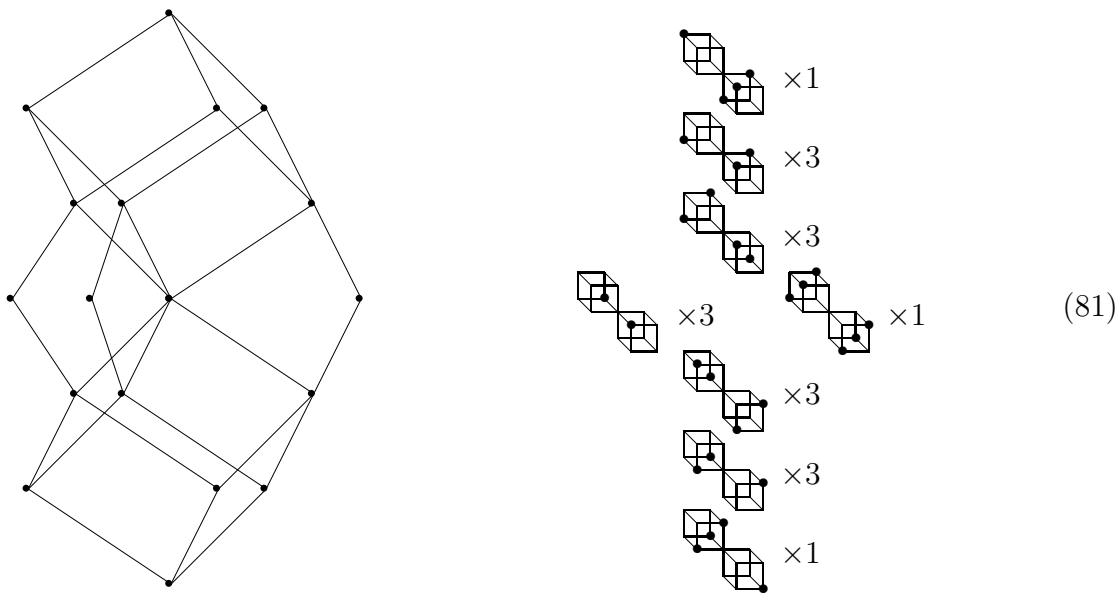
\begin{equation}\begin{minipage}[c]{35pc}\label{d.3cx_outline}
\begin{picture}(250,225)

\put(220,0){
\put(60,0){\twocube{
\put(30,0){\bigdot}
\put(5,15){\bigdot}
\put(10,20){\bigdot}
\put(15,25){\bigdot}
}\put(5,10){$\times 1$}}

\put(60,30){\twocube{
\put(30,10){\bigdot}
\put(5,15){\bigdot}
\put(10,20){\bigdot}
}\put(5,10){$\times 3$}}

\put(60,60){\twocube{
\put(20,0){\bigdot}
\put(30,10){\bigdot}
\put(10,20){\bigdot}
\put(5,25){\bigdot}
}\put(5,10){$\times 3$}}

\put(20,90){\twocube{
\put(20,10){\bigdot}
\put(10,20){\bigdot}
}\put(5,10){$\times 3$}}

\put(100,90){\twocube{
\put(20,0){\bigdot}
\put(25,5){\bigdot}
\put(30,10){\bigdot}
\put(0,20){\bigdot}
\put(5,25){\bigdot}
\put(10,30){\bigdot}
}\put(5,10){$\times 1$}}

\put(60,120){\twocube{
\put(20,10){\bigdot}
\put(25,5){\bigdot}
\put(0,20){\bigdot}
\put(10,30){\bigdot}
}\put(5,10){$\times 3$}}

\put(60,150){\twocube{
\put(20,10){\bigdot}
\put(25,15){\bigdot}
\put(0,20){\bigdot}
}\put(5,10){$\times 3$}}

\put(60,180){\twocube{
\put(15,5){\bigdot}
\put(20,10){\bigdot}
\put(25,15){\bigdot}
\put(0,30){\bigdot}
}\put(5,10){$\times 1$}}
}

\put(90,110){
\multiput(0,0)(-36,36){2}{
\multiput(0,0)(54,36){2}{\put(0,0){\line(-1,2){18}}
	\multiput(0,0)(-18,36){2}{\middot}}
\multiput(0,0)(-18,36){2}{\put(0,0){\line(3,2){54}}}
}
\multiput(0,0)(54,36){2}{\multiput(0,0)(-18,36){2}{
\put(0,0){\line(-1,1){36}}}}

\put(-36,36){\line(-2,-3){24}}\put(-60,0){\middot}
\put(-18,36){\line(-1,-3){12}}\put(-30,0){\middot}
\put(54,36){\line(1,-2){18}}\put(72,0){\middot}

\multiput(0,0)(-36,-36){2}{
\multiput(0,0)(54,-36){2}{\put(0,0){\line(-1,-2){18}}
	\multiput(0,0)(-18,-36){2}{\middot}}
\multiput(0,0)(-18,-36){2}{\put(0,0){\line(3,-2){54}}}
}
\multiput(0,0)(54,-36){2}{\multiput(0,0)(-18,-36){2}{
\put(0,0){\line(-1,-1){36}}}}

\put(-36,-36){\line(-2,3){24}}
\put(-18,-36){\line(-1,3){12}}
\put(54,-36){\line(1,2){18}}
}
\end{picture}
\vspace{1em}
\end{minipage}\end{equation}
The reader familiar with lattice theory will recognize the diagram
at the left as the free distributive lattice on three generators
\cite[Figure~19, p.\,84]{GLT3}.
This is because the distributive lattice of proper nonempty down-sets of
the set of vertices of our cube, under union and intersection,
is freely generated by the three down-sets
\raisebox{-.5em}{\begin{picture}(25,15)
\multiput(10,0)(-5,5){2}{
\multiput(0,0)(0,10){2}{\line(1,0){10}\bigdot}
\multiput(0,0)(10,0){2}{\line(0,1){10}}}
\multiput(10,0)(10,0){2}{\multiput(0,0)(0,10){2}{\shortdiag}}
\end{picture}},
\raisebox{-.5em}{\begin{picture}(25,15)
\multiput(10,0)(-5,5){2}{
\multiput(0,0)(0,10){2}{\line(1,0){10}}
\multiput(0,0)(10,0){2}{\line(0,1){10}\bigdot}}
\multiput(10,0)(10,0){2}{\multiput(0,0)(0,10){2}{\shortdiag}}
\end{picture}}
and
\raisebox{-.5em}{\begin{picture}(25,15)
\multiput(10,0)(-5,5){2}{
\multiput(0,0)(0,10){2}{\line(1,0){10}}
\multiput(0,0)(10,0){2}{\line(0,1){10}}}
\multiput(10,0)(10,0){2}{\multiput(0,0)(0,10){2}
	{\shortdiag\put(0,0){\bigdot}}}
\end{picture}}.
(These free generators are the three ``outer'' elements
at the middle level of the diagram on the left above,
corresponding, in the diagram on the right, to the
picture at that level marked~$\times 3.$
As constructions on our triple complex, they represent
the classical homology constructions corresponding to the
three axial directions in that complex.)
\vspace{.3em}

I have not investigated what extramural maps and exact sequences relate
these $18$ constructions.
I will not propose iconic notations for them like those we have used
in studying double complexes;
probably the best notation would, rather, involve indexing them by
expressions for the elements of a free distributive lattice;
e.g., $h_{(x_1\wedge x_2)\vee x_3}(A_{ijk})$ or
$h((x\wedge y)\vee z;\,A_{ijk}),$
where $x_1,$ $x_2,$ $x_3,$ or $x,$ $y,$ $z,$ denote the free generators.
For informal purposes, though, something like $h(
\raisebox{-.5em}{\begin{picture}(25,15)
\multiput(8,0)(-5,5){2}{
\multiput(0,0)(0,10){2}{\line(1,0){10}}
\multiput(0,0)(10,0){2}{\line(0,1){10}}}
\multiput(8,0)(10,0){2}{\multiput(0,0)(0,10){2}{\shortdiag}}
\put(8,10){\bigdot}
\put(13,5){\bigdot}
\end{picture}}\kern-.4em,\,A_{ijk})$\vspace{.3em}
might occasionally be convenient (as long as we don't go beyond
triple complexes).

Are these objects likely to be of use?
I don't know!

\subsection{Kernel and image ratios}\label{SS.ratios}
J.\,Lambek \cite{JL} (cf.\ \cite[Lemma~III.3.1]{H+S}) associates to any
commuting square,
\begin{equation}\begin{minipage}[c]{30pc}\label{d.comm_sq}
\begin{picture}(50,65)
\put(5,50){$P$}
\put(18,53){\vector(1,0){25}\put(-17,3){$a$}}
\put(18,45){\vector(1,-1){25}\put(-12,-12){$f$}}
\put(45,50){$R$}
\put(10,45){\vector(0,-1){25}\put(-8,-13){$b$}}
\put(50,45){\vector(0,-1){25}\put(1,-13){$c$}}
\put(5,10){$Q$}
\put(18,13){\vector(1,0){25}\put(-17,-10){$d$}}
\put(45,10){$S$}
\put(100,30){where $ca\ =\ f\ =\ db$}
\end{picture}
\end{minipage}\end{equation}
two objects, which he calls the \emph{kernel ratio}
and the \emph{image ratio} of the square.
To bring out the similarity to the concepts of this
paper, let me name them
\begin{equation}\begin{minipage}[c]{30pc}\label{d.P**S}
$P_*\ =\ \r{Ker}\,f/(\r{Ker}\,a\,+\,\r{Ker}\,b),\\[.3em]
^*S\ =\ (\r{Im}\,c\ \cap\,\r{Im}\ d)\,/\,\r{Im}\,f\,.$
\end{minipage}\end{equation}
In fact, if the given square is embedded in a double complex
which is vertically \emph{and} horizontally exact at $P,$
respectively at $S,$ then we see that $P_*=P_\bx\,,$ respectively,
$^* S={^\bx S}.$

Now suppose we have a commuting diagram with exact rows,
\begin{equation}\begin{minipage}[c]{30pc}\label{d.2comm_sq}
\begin{picture}(50,65)
\put(5,50){$P$}
\put(18,53){\vector(1,0){25}}
\put(45,50){$R$}
\put(58,53){\vector(1,0){25}}
\put(85,50){$T$}
\put(10,45){\vector(0,-1){25}}
\put(50,45){\vector(0,-1){25}\put(1,-13){$c$}}
\put(90,45){\vector(0,-1){25}}
\put(5,10){$Q$}
\put(18,13){\vector(1,0){25}}
\put(45,10){$S$}
\put(58,13){\vector(1,0){25}}
\put(85,10){$U\,.$}
\end{picture}
\end{minipage}\end{equation}
Then we can extend this diagram, by putting in the kernel
and cokernel of $c,$ to a double complex exact in both
directions at $R$ and $S.$
Hence Corollary~\ref{C.extra-iso}, applied to $c,$ gives us
\begin{equation}\begin{minipage}[c]{30pc}\label{d.*+bx}
$R_*\ \cong\ R_\bx\ \cong\ ^\bx S\ \cong\ ^* S.$
\end{minipage}\end{equation}
The isomorphism $R_*\cong{^* S}$ is proved by Lambek \cite{JL}
and used to get other results, as I use Corollary~\ref{C.extra-iso} in
\S\ref{S.easy} above.
The constructions $(\ )_*$ and $^* (\ )$ have the advantage
of being definable with reference to a smaller diagram
than my $(\ )_\bx$ and $^\bx (\ ).$
They share with these constructions the property of vanishing on any
\emph{doubly exact} double complex with finite support.
But one doesn't seem to be able to do anything with them
without some exactness assumptions.
With such assumptions, one gets extramural \emph{iso}morphisms
as in Corollary~\ref{C.extra-iso}, but
without them, one does not have analogs of the
extramural \emph{homo}morphisms of Lemma~\ref{L.salamander}.
\subsection{Non-abelian groups}\label{SS.nonabelian}
Lambek proves the results referred to above for not necessarily
abelian groups, through he applies them in abelian situations.
Note that for non-abelian groups, the exactness of the top row
of~(\ref{d.2comm_sq}), or something
similar, is needed to conclude that $^* S$ will
be a group, i.e., that the denominator in the definition will
be a \emph{normal} subgroup of the numerator.
Without that  condition, $^* S$ is a ``homogeneous space''.
I likewise noticed when first working out the Salamander Lemma
that a version could be stated for not-necessarily-abelian
groups, but there were even worse difficulties -- e.g.,
in Definition~\ref{D.4things}, $A_\bx$ would just be a pointed set,
the quotient of the group $\r{Ker}\,q$ by the left action of
the subgroup $\r{Im}\,c$ and the right action of
the subgroup $\r{Im}\,d.$
However, it would certainly be nice to have a tool like the
Salamander Lemma for proving noncommutative versions of
basic diagram-chasing lemmas, when these hold.

Leicht \cite{JBL}, Kopylov \cite{kopylov}, and others have
given more general conditions on
a category under which Lambek's result holds.
In \cite{JL.b+s}, Lambek gets related results for varieties
of algebras in the sense of universal
algebra satisfying an appropriate Mal'cev-type condition.

\section{Two exercises}\label{S.exercise}

I have left many key calculations in this note
to the reader, including the verification of
the Salamander Lemma (Lemma~\ref{L.salamander}) itself.
My ``sketches of proofs'' can likewise
be regarded ``exercises with hints''.
Here are two further interesting exercises.

\subsection{Building up finite exact double complexes}\label{S.build}
Let $\A$ be an abelian category, and let $\A^\#$ denote the
abelian category of double complexes in $\A,$ where a morphism
$f:A\to B$ is a family of morphisms $f_{i,r}:A_{i,r}\to B_{i,r}$
commuting with $\delta_1$ and $\delta_2.$
(We are not admitting subscript-shifting morphisms here.)
Let $\r{FX}\subseteq\r{Ob}(\A^\#)$ (standing for ``finite exact'')
be the class of double complexes with
only finitely many nonzero objects, and all rows and columns exact,
and $\r{EX}\subseteq\r{FX}$ (``elementary exact'')
be the class of double complexes of the form
\begin{equation}\begin{minipage}[c]{15pc}\label{d.4A}
\begin{picture}(120,135)
\put(40,5){$0$}
\put(80,5){$0$}
\put(0,45){$0$}
\put(38,45){$A$}
\put(78,45){$A$}
\put(120,45){$0$}
\put(0,85){$0$}
\put(38,85){$A$}
\put(78,85){$A$}
\put(120,85){$0$}
\put(40,125){$0$}
\put(80,125){$0$}
\multiput(11,49)(0,40){2}{\multiput(0,0)(40,0){3}{\vector(1,0){24}}}
\multiput(44,121)(40,0){2}{\multiput(0,0)(0,-40){3}{\vector(0,-1){24}}}
\end{picture}
\end{minipage}\end{equation}
placed at arbitrary locations in the grid,
where the maps among all copies of $A$ are the identity.\\[.5em]
(a)  Show that $\r{FX}$ is the least subclass of $\r{Ob}(\A^\#)$
containing $\r{EX}$ and closed therein under extensions.

Hint:  Given an object of $\r{FX},$ prove that one can map it
epimorphically to an object of $\r{EX}$ as suggested below.
\begin{equation}\begin{minipage}[c]{33pc}\label{d.onto_4D}
\begin{picture}(120,140)
\put(80,5){$0$}
\put(38,45){$C$}
\put(78,45){$D$}
\put(120,45){$0$}
\put(38,85){$A$}
\put(78,85){$B$}
\multiput(11,49)(0,40){2}{\multiput(0,0)(40,0){3}{\vector(1,0){24}}}
\multiput(44,121)(40,0){2}{\multiput(0,0)(0,-40){3}{\vector(0,-1){24}}}

\put(171,65){\thicklines\vector(1,0){30}\vector(1,0){4}}

\put(240,0){
\put(40,5){$0$}
\put(80,5){$0$}
\put(0,45){$0$}
\put(38,45){$D$}
\put(78,45){$D$}
\put(120,45){$0$}
\put(0,85){$0$}
\put(38,85){$D$}
\put(78,85){$D$}
\put(120,85){$0$}
\put(40,125){$0$}
\put(80,125){$0$}
\multiput(11,49)(0,40){2}{\multiput(0,0)(40,0){3}{\vector(1,0){24}}}
\multiput(44,121)(40,0){2}{\multiput(0,0)(0,-40){3}{\vector(0,-1){24}}}
\put(120,10){.}
}
\end{picture}
\end{minipage}\end{equation}

To get epimorphicity at $A,$ show that $^\bx D=0.$\\[.5em]
(b)  In contrast, show by example that for a double complex which
is \emph{not} assumed to be finite, but nonetheless has
exact rows and columns, and has an object
$D$ with zero objects immediately below it and to its right,
the map of~(\ref{d.onto_4D}) may fail to be an epimorphism.
\vspace{.5em}

I expect that an epimorphicity result analogous
to~(\ref{d.onto_4D}) should hold for finite exact triple (and
higher) complexes.
If so, one could get the analog of~(a) above, and
deduce from this that the $18$ constructions of \S\ref{SS.3-cx}
(corresponding to the diagrams of~(\ref{d.3cx_outline}))
all give zero at objects of a finite triple complex exact in
all three directions.

\subsection{Complexes with a twist}\label{S.twist}
(a)  Suppose we are given a commutative diagram
\begin{equation}\begin{minipage}[c]{35pc}\label{d.twist}
\begin{picture}(350,230)
\put(11,11){
\multiput(40,40)(160,-40){2}{\multiput(0,0)(0,160){2}{
	\multiput(0,0)(40,0){3}{$0$}}}
\multiput(14,83)(0,40){2}{\multiput(0,0)(40,0){8}{\vector(1,0){24}}}
\multiput(14,163)(200,-120){2}{\multiput(0,0)(40,0){3}{%
	\vector(1,0){24}}}
\multiput(45,75)(160,-40){2}{\multiput(0,0)(40,0){3}{
	\multiput(0,0)(0,40){4}{\vector(0,-1){24}}}}

\qbezier(134,163)(165,163)(165,133)\put(165,133){\vector(0,-1){2}}
\put(165,115){\vector(0,-1){24}}
\qbezier(165,75)(165,43)(195,43)\put(197,43){\vector(1,0){2}}
}

\put(10,11){
\put(40,160){$A$}
\put(80,160){$B$}
\put(120,160){$C$}

\put(40,120){$L$}
\put(79,120){$M$}
\put(120,120){$N$}
\put(160,120){$E$}
\put(199,120){$P$}
\put(240,120){$Q$}
\put(279,120){$R$}

\put(40,80){$U$}
\put(80,80){$V$}
\put(119,80){$W$}
\put(160,80){$F$}
\put(199,80){$X$}
\put(240,80){$Y$}
\put(279,80){$Z$}

\put(199,40){$H$}
\put(240,40){$J$}
\put(279,40){$K$}
}
\end{picture}\end{minipage}\end{equation}
in which the columns beginning and ending with $0$ are short
exact sequences, while the three ``rows'' (the two ordinary ones, and
the one that makes a detour from the top
to the bottom) are merely assumed to be complexes.

Obtain a long exact sequence of (mostly horizontal) homology objects
\begin{equation}\begin{minipage}[c]{34pc}\label{d.mostly_hor}
$\cdots{\to}
V\kern-.10em\hor\,{\to}
C\kern-.10em\hor\,{\to}
N\kern-.10em\hor\,{\to}
W\kern-.20em\hor\,{\to}
E\ver\,{\to}
E\kern-.1em\hor\,{\to}
F\kern-.15em\hor\,{\to}
F\ver\,{\to}
P\kern-.15em\hor\,{\to}
X\kern-.20em\hor\,{\to}
H\kern-.20em\hor\,{\to}
Q\hor\,{\to}\cdots.$
\end{minipage}\end{equation}

(Suggestion: Square off the curved arrows
in~(\ref{d.twist}) by inserting the
kernel $D$ and the cokernel $G$ of the arrow $E\to F,$
and add zeroes above, respectively, below these,
getting an ordinary double complex with all columns exact.
Apply to this the idea used in deriving Lemma~\ref{L.ijk}.)\\[.5em]
(b)  Consider any double complex~(\ref{d.some_labeled}),
in which we have labeled some of
the objects to the upper left and lower right of one
arrow, $F\to G$:

\begin{equation}\begin{minipage}[c]{20pc}\label{d.some_labeled}
\begin{picture}(240,180)
\put(036,033){\middot}
\put(066,033){\middot}
\put(096,033){\middot}
\put(120,030){$M$}
\put(156,033){\middot}
\put(186,033){\middot}

\put(036,063){\middot}
\put(066,063){\middot}
\put(096,063){\middot}
\put(120,060){$K$}
\put(150,060){$L$}
\put(186,063){\middot}

\put(030,090){$D$}
\put(060,090){$E$}
\put(090,090){$F$}
\put(120,090){$G$}
\put(150,090){$H$}
\put(180,090){$J$}

\put(036,123){\middot}
\put(060,120){$B$}
\put(090,120){$C$}
\put(126,123){\middot}
\put(156,123){\middot}
\put(186,123){\middot}

\put(036,153){\middot}
\put(066,153){\middot}
\put(090,150){$A$}
\put(126,153){\middot}
\put(156,153){\middot}
\put(186,153){\middot}

\multiput(15,34)(0,30){5}{\multiput(0,0)(30,0){7}{\vector(1,0){12}}}
\multiput(36,24)(30,0){6}{\multiput(0,0)(0,30){6}{\vector(0,-1){12}}}
\put(210,10){.}
\end{picture}
\end{minipage}\end{equation}

Obtain from this a diagram
\begin{equation}\begin{minipage}[c]{35pc}\label{d.twist+oplus}
\begin{picture}(350,230)
\put(20,10){
\multiput(40,40)(190,-40){2}{\multiput(0,0)(0,160){2}{
	\multiput(0,0)(70,0){2}{$0$}}}
\multiput(45,75)(190,-40){2}{\multiput(0,0)(70,0){2}{
	\multiput(0,0)(0,40){4}{\vector(0,-1){24}}}}

\qbezier(124,163)(175,163)(175,133)\put(175,133){\vector(0,-1){2}}
\put(175,115){\vector(0,-1){24}}
\qbezier(175,75)(175,42)(226,43)\put(226,43){\vector(1,0){2}}

\put(-7,163){\vector(1,0){44}}
\put(39,160){$D$}
\put(53,163){\vector(1,0){54}}
\put(109,160){$E$}

\put(-7,123){\vector(1,0){22}}
\put(15,120){$A\oplus B\oplus D$}
\put(78,123){\vector(1,0){14}}
\put(96,120){$C\oplus E$}
\put(133,123){\vector(1,0){34}}
\put(169,120){$F$}
\put(183,123){\vector(1,0){44}}
\put(229,120){$K$}
\put(243,123){\vector(1,0){44}}
\put(288,120){$L\oplus M$}
\put(325,123){\vector(1,0){42}}

\put(-7,83){\vector(1,0){30}}
\put(26,80){$A\oplus B$}
\put(63,83){\vector(1,0){40}}
\put(109,80){$C$}
\put(123,83){\vector(1,0){44}}
\put(169,80){$G$}
\put(183,83){\vector(1,0){29}}
\put(215,80){$H\oplus K$}
\put(255,83){\vector(1,0){19}}
\put(278,80){$J\oplus L\oplus M$}
\put(341,83){\vector(1,0){24}}

\put(229,40){$H$}
\put(243,43){\vector(1,0){54}}
\put(299,40){$J$}
\put(313,43){\vector(1,0){54}}
}
\end{picture}
\end{minipage}\end{equation}
taking the arrows to be the sums of ``all available'' morphisms.
(E.g., the arrow $E\to C\oplus E$ is just the inclusion, but the
arrow $C\oplus E\to F$ is the sum of the given
arrows $C\to F$ and $E\to F.)$
Verify that, after a bit of sign-tweaking,~(\ref{d.twist+oplus})
satisfies the hypotheses of part~(a).
Write out, for this diagram, the middle six terms of the exact
sequence~(\ref{d.mostly_hor}), then one more on either side.

Thus we see that we can extend the $\!6\!$-term ``salamanders''
of Lemma~\ref{L.salamander}
to longer exact sequences, if we are willing to define more
complicated auxiliary objects.\\[.5em]
(c) Do these new constructions still have the property of being zero
on \emph{bounded exact} double complexes?

\section{Acknowledgements, and a final remark}\label{S.ackn}
My parents, Lester and Sylvia Bergman, worked together in scientific
photography and illustration; and though they never instructed me
in the latter, I acquired from them my love of an
effective visual display.
(But the responsibility for the often inconsistent
\TeX\ coding underlying the diagrams of this note is my own.)

Most of this work was done, and a rough draft
written, in the Fall of 1972, when I
was supported by an Alfred P.\ Sloan Research Fellowship, and
was a guest at the University of Leeds' stimulating Ring Theory Year.
In 2007 Anton Geraschenko digitized, and, with my permission,
put online, a copy of that draft.
His enthusiasm for the material helped spur me to create this
better version.
I am also indebted to the referee for several thoughtful suggestions,
though shortage of time has prevented me from following many of them.

I have not attempted to make this note a ``definitive''
development of the subject: that should be left to those
who actually work with double complexes, and can judge
how best to develop the material.
Hence the variation, in the present note,
between sections where objects are denoted
by arbitrary letters and those where they are distinguished by double
subscripts (as seemed to give the most readable
presentation of one or another topic), the lack of
any general notation for intramural and extramural maps, and the
choice of vertical arrows that point downward (as in
most familiar lemmas proved by diagram chasing) rather than
upward (as might be preferable based on systematic considerations).

% - - - - - - - - - - - - - - - - - - - - - - - - - - - - - -

\end{document}